\documentclass[12pt, twoside, leqno]{article}
\usepackage{amsmath,amsthm,amscd}
\usepackage{amssymb}
\usepackage{enumerate}
\usepackage{amsthm}
\newtheorem{thm}{Theorem}[section]
\newtheorem{cor}[thm]{Corollary}
\newtheorem{lem}[thm]{Lemma}
\newtheorem{prop}[thm]{Proposition}
\newtheorem{rem}[thm]{Remark}
\theoremstyle{definition}

\newcommand{\Real}{{\mathbb{R}}}
\newcommand{\Com}{{\mathbb{C}}}
\newcommand{\Zet}{{\mathbb{Z}}}
\newcommand{\Nat}{{\mathbb{N}}}

\newcommand{\field}{{\mathbb{K}}}
\newcommand{\Pro}{{\mathbb{P}}}
\newcommand{\Oc}{{\mathcal{O}}}

\newcommand{\Li}{{\mathcal{L}}}

\newcommand{\Id}{{\mathcal{I}}}
\newcommand{\Gl}{{\mathrm{GL}}}

\newcommand{\Z}{{\mathcal{Z}}}

\newcommand{\U}{{\mathcal{U}}}

\newcommand{\Lie}{{\mathrm{Lie}}}
\newcommand{\ord}{{\mathrm{ord}}}
\newcommand{\Pic}{{\mathrm{Pic}}}
\newcommand{\G}{{\overline{G}}}
\newcommand{\Hb}{{\overline{H}}}
\newcommand{\Lb}{{\overline{L}}}

\newcommand{\Spec}{{\mathrm{Spec}}}

\newcommand{\gfr}{{\mathfrak{g}}}
\newcommand{\bfr}{{\mathfrak{b}}}
\newcommand{\Ac}{{\mathcal{A}}}
\newcommand{\Dc}{{\mathcal{D}}}
\newcommand{\Bc}{{\mathcal{B}}}

\newcommand{\Af}{{\mathbb{A}}}
\newcommand{\x}{{\mathbf{x}}}
\newcommand{\y}{{\mathbf{y}}}
\newcommand{\Ib}{{\mathbf{In}}}
\newcommand{\Ih}{{^hI}}
\newcommand{\Bb}{{\mathbf{B}}}
\newcommand{\Tb}{{\mathbf{T}}}
\newcommand{\Db}{{\mathbf{D}}}
\newcommand{\Cb}{{\mathbf{C}}}
\newcommand{\Eb}{{\mathbf{E}}}
\newcommand{\Der}{{\mathrm{Der}}}
\newcommand{\Ad}{{\mathrm{Ad}}}
\numberwithin{equation}{section}

\begin{document}
\setlength{\parindent}{0pt}
\baselineskip=17pt

\title{On the multiplicity estimates}

\author{Mario  Huicochea\\
Facultad de Ciencias, UNAM\\ 
E-mail: dym@cimat.mx
}

\date{}

\maketitle
\begin{abstract}
 In this paper we show some multiplicity estimates theorems for a connected algebraic group (not necessarily commutative) $G$  over an algebraically closed subfield of $\Com$. More specifically, under particular assumptions on the parameters and the points where the polynomial has high order with respect to a Lie subalgebra of the Lie algebra associated to $G$, we present a series of results where we find   obstruction varieties  with different properties. Some of the results obtained in this paper improve the multiplicity estimates theorem for arbitrary connected algebraic groups that already exist, see \cite [Thm. 0.3]{Nakamaye1}. 
\end{abstract}
\section{Introduction}
In transcendence  theory, the multiplicity estimates theorems have been quite important.  Some of the main results of transcendence theory have used multiplicity estimates theorems in fundamental parts of their proofs, see for instance \cite[Hauptsatz]{Wustholz2}, \cite[Thm. 4.1]{Waldschmidt}  and \cite[Thm. 3]{Gaudron}.  One of the first multiplicity estimates theorem was obtained by Nesterenko \cite{Nesterenko}. In the next few years, several improvements were done, see for example \cite{BrownawellMasser}, \cite{Masser-Wustholz1} or \cite{Masser-Wustholz2}. Some years later W\"{u}stholz \cite[Main Thm.]{Wustholz1} and Philippon \cite[Thm. 2.1]{Philippon} published  two breakthrough results; these papers were quite important in the developments of number theory  since they have had several applications in transcendence theory. Also  W\"{u}stholz and Philippon results were improved in different directions  in several papers, see for instance \cite[Thm 0.3]{Nakamaye1}, \cite[Thm. 1.1]{Wustholz3}, \cite[Thm. 1]{Nakamaye2} and \cite[Thm. 1]{Fischler Nakamaye}. In particular we shall be interested in the point of view of Nakamaye in \cite[Thm. 0.3]{Nakamaye1}. The results of  W\"{u}stholz and Philippon were done for connected commutative algebraic groups and Nakamaye remarked  that most of the tools, jointly with some technical assumptions,  used in the proof of \cite[Thm. 0.3]{Nakamaye1} are  
generalizable to connected algebraic groups. The  goal of this paper is to continue with the study of multiplicity estimates theorems for noncommutative algebraic groups.

Let $\field$ be an algebraically closed subfield of $\Com$ and  $G$  a connected algebraic group over $\field$ of dimension $n$. In Section 2 we construct  a $G-$biequivariant compactification $\G$ and  a closed embedding $\phi:\G\rightarrow \Pro^N$. We denote by $I(\G)$ the set of polynomials of $\field[x_0,\ldots,x_N]$ which vanish in  $\phi(\G)$.  We will denote by  $\gfr$ the  corresponding $\field-$Lie algebra of $G$; in other words $\gfr$ is the set of left invariant elements of  $\mathrm{Der}(\Oc_G,\Oc_G)$ and it may be identified with the left invariant elements of  $\mathrm{Der}(\Oc_\G,\Oc_\G)$, see Section 3.  When $I$ is a homogeneous ideal of $\field[x_0,\ldots,x_N]$, $\Z(I)\subseteq \Pro^N$ denotes its zero set.  For all $k\in\{0,\ldots,N\}$, set $G_k:=\{z\in\G:\; x_k(\phi(z))\neq 0\}$ and assume without loss of generality that $1\in \bigcap_{k=0}^nG_k$.  Let $\bfr$ be  a Lie subalgebra of $\gfr$ and $\{\Delta_1, \ldots, \Delta_d\}$  a fixed basis of $\bfr$; denote by $\U(\bfr)$ the universal enveloping algebra of $\bfr$. We denote by $\U(\bfr,T)$ the $\field-$subspace of $\U(\bfr)$ generated by   $\Delta^{t_1}_1\ldots \Delta^{t_d}_d$ where $t_1,\ldots, t_d\in\Nat\cup\{0\}$ and $\sum_{i=1}^d t_i\leq T$. For any element $\Delta:=\Delta^{t_1}_1\ldots \Delta^{t_d}_d$,
 define
\begin{equation*} 
  \Delta(G_k):=\overbrace{\Delta_{d}(G_k)\circ\ldots\circ \Delta_{d}(G_k)}^{t_d}\circ\ldots\circ\overbrace{\Delta_{1}(G_k)\circ\ldots\circ \Delta_{1}(G_k)}^{t_1};
  \end{equation*} 
 then the definition of $\Delta(G_k)$ for arbitrary $\Delta\in\U(\bfr)$ is extended by linearity. Given  $g\in  G$ let $Q_0,\ldots,Q_N\in\field[x_0,\ldots,x_N]$ be homogeneous polynomials of the same degree and $U\subseteq G$  a neighbourhood  of $1$  such that $\phi(gz)=\big[Q_0(z):\ldots:Q_N(z)\big]$ for all  $z\in U$  and write $Q:=(Q_0,\ldots, Q_N)$. For $P\in\field[x_0,\ldots, x_N]\setminus I(\G)$ homogeneous, define 
    $\ord_{g}(\bfr,P)$ as the minimum $T\in\Nat\cup\{0\}$  such that there exists $\Delta\in\U(\bfr,T)$ satisfying
\begin{equation*}    
 \Delta(G_0)\bigg(P\circ Q\bigg(\frac{x_0}{x_0},\ldots,\frac{x_N}{x_0}\bigg)\bigg)\Bigg|_{[x_0:\ldots:x_N]=\phi(1)}\neq 0;
 \end{equation*}   the definition of $\ord_{g}(\bfr,P)$ depends neither on $Q_0,\ldots, Q_N$ nor on $U$. We identify $G(\field):=\mathrm{Hom}_{\Spec(\field)}(\Spec(\field),G)$ with the closed points of $G$ and we consider it a subset of $G(\Com):=\mathrm{Hom}_{\Spec(\field)}(\Spec(\Com),G)$. $G(\Com)$ has a $\Com-$Lie group structure with $\Lie(G(\Com))\cong\gfr\otimes_\field\Com$. Call $B$ the connected analytic subgroup of $G(\Com)$ corresponding to the Lie subalgebra $\bfr\otimes_\field\Com$ of $\Lie(G(\Com))$. For an irreducible  subvariety  $W$ of $\G$ and $w\in W(\Com)$ such that $W(\Com)\cap Bw$ is transverse at $w$, call 
\begin{equation*} 
 \tau(W):=\dim(B)-\dim(W(\Com)\cap Bw). 
\end{equation*} 
 If $V$ is an irreducible variety of  $\phi(\G)$, set $\tau(V):=\tau(\phi^{-1}(V))$.  For a projective variety $V$ embedded in $\Pro^N$, we denote by $\deg(V)$ the degree of $V$. For a finite subset $\Sigma_1$ of $G$ containing $1$, define  $\Sigma_0=\{1\}$ and $\Sigma_S:=\big\{\prod_{k=1}^S g_k:\;g_1,\ldots,g_S\in \Sigma_1\big\}$ for $S\in\Nat$.  Finally, for $x\in\Real$ we denote by $[x]$  the largest integer less than or equal to  $x$. The first main result of this paper is the following.
\begin{thm}
\label{R1}
Let   $P\in\field[x_0,\ldots,x_N]\setminus I(\G)$ be homogeneous of degree $D$. Assume  that $\ord_{g}(\bfr,P)\geq T+1$ for all $g\in \Sigma_S$ and  $D\geq\sum_{i=0}^{S}(|\Sigma_1|-1)^i$. Then there exists $c_1$ independent of $S,T,D,\bfr, \Sigma_1$ and $P$ with the following property: for all  $ d_0\in\{1,\ldots, n\}$ there is an irreducible subvariety $W$ of $\G$ such that
\begin{enumerate}
\item[i)]$\dim(W)\leq d_0$.
\item[ii)]$1\in W\cap G$.
\item[iii)]$\phi(W)\subseteq \Z(P)$.
\item[iv)]If $N_W$ is the number of different cosets $gW$ for all $g\in \Sigma_{[\frac{S}{n}]}$, then 
\begin{equation*}
N_W\binom{\big[\frac{T}{n}\big]+\tau(W)}{\tau(W)}\deg(\phi(W))\leq c_1D^{n-\dim(W)}.
\end{equation*}
\end{enumerate}
\end{thm}
If we remove the hypothesis $D\geq\sum_{i=0}^{S}(|\Sigma_1|-1)^i$ from Theorem \ref{R1}, we cannot assure, given a $ d_0\in\{1,\ldots, n\}$, the existence of $W$ satisfying i). Nonetheless, also in this case we may find $W$ satisfying ii)-iv); in other words we shall demonstrate the following statement. 
\begin{thm}
\label{R2}
Let   $P\in\field[x_0,\ldots,x_N]\setminus I(\G)$ be homogeneous of degree $D$. Assume  that $\ord_{g}(\bfr,P)\geq T+1$ for all $g\in\Sigma_S$. Then there exists $c_2$ independent of $S,T,D,\bfr,\Sigma_1$ and $P$ with the following property: there is an irreducible subvariety $W$ of $\G$ such that
\begin{enumerate}
\item[i)]$1\in W\cap G$.
\item[ii)]$\phi(W)\subseteq \Z(P)$.
\item[iii)]If $N_W$ is the number of different cosets $gW$ for all $g\in\Sigma_{[\frac{S}{n}]}$, then 
\begin{equation*}
N_W\binom{\big[\frac{T}{n}\big]+\tau(W)}{\tau(W)}\deg(\phi(W))\leq c_2D^{n-\dim(W)}.
\end{equation*}
\end{enumerate}
\end{thm}
In the previous theorems, we find irreducible subvarieties of $\G$ satisfying some properties; it is natural to ask whether we can assure that $W$ is in an interesting family of varieties. In the next results, we shall show that under certain  hypothesis $W$ may be assumed to be the closure of an  irreducible normal  algebraic subgroup in $\G$. Moreover, the following theorem could have  application in transcendence theory for noncommutative algebraic groups.   
\begin{thm}
\label{R3}
Let   $P\in\field[x_0,\ldots,x_N]\setminus I(\G)$ be homogeneous of degree $D$. Assume  that $\ord_{g}(\bfr,P)\geq T+1$ for all $g\in \Sigma_{S}$,   $gB=Bg$ for all $g\in \Sigma_1$ and  $\sum_{i=0}^{S}(|\Sigma_1|-1)^i\leq D$. Then there exists $c_3$ independent of $S,T,D,\bfr,\Sigma_1$ and $P$ with the following property: for all  $ d_0\in\{1,\ldots, n\}$ there is an irreducible algebraic subgroup $H$ of $G$ such that
\begin{enumerate}
\item[i)]$\dim(H)\leq d_0$.
\item[ii)]$\phi(H)\subseteq \Z(P)$.
\item[iii)]$H$ is a normal subgroup of $G$.
\item[iv)]Let $\Hb$ be the closure of $H$ in $\G$. If $N_H$ is the number of different cosets $Hg$ for all $g\in\Sigma_{[\frac{S}{n}]}$, then 
\begin{equation*}
N_H\binom{\big[\frac{T}{n}\big]+\tau(\Hb)}{\tau(\Hb)}\deg(\phi(\Hb))\leq c_3D^{n-\dim(H)}.
\end{equation*}
\end{enumerate}
\end{thm}
The last main result will be a mix between Theorem \ref{R2} and Theorem \ref{R3}; concretely, we remove the assumption $D\geq \sum_{i=0}^{S}(|\Sigma_1|-1)^i$ from Theorem \ref{R3}, nevertheless we  assure that  the obstruction variety  is the closure of a normal algebraic subgroup. This theorem generalizes  \cite[Thm. 0.3]{Nakamaye1}; specifically,  it is assumed in \cite[Thm. 0.3]{Nakamaye1} that $B$ is the image of $\Com-$Lie group morphism $\psi:\Com^d\rightarrow G(\Com)$ and therefore this image is a commutative subgroup of $G(\Com)$, see \cite[p. 157]{Nakamaye1}. In the following theorem we show that this assumption is not needed.
\begin{thm}
\label{R4}
Let   $P\in\field[x_0,\ldots,x_N]\setminus I(\G)$ be homogeneous of degree $D$. Assume  that $\ord_{g}(\bfr,P)\geq T+1$ for all $g\in\Sigma_S$ and $gB=Bg$ for all $g\in \Sigma_1$. Then there exists $c_4>0$ independent of $S,T,D,\bfr,\Sigma_1$ and $P$  with the following property:  there is an irreducible algebraic subgroup $H$ of $G$ such that
\begin{enumerate}
\item[i)]$\phi(H)\subseteq \Z(P)$.
\item[ii)]$H$ is a normal subgroup of $G$.
\item[iii)]Let $\Hb$ be the closure of $H$ in $\G$. If $N_H$ is the number of different cosets $Hg$ for all $g\in\Sigma_{[\frac{S}{n}]}$, then 
\begin{equation*}
N_H\binom{\big[\frac{T}{n}\big]+\tau(\Hb)}{\tau(\Hb)}\deg(\phi(\Hb))\leq c_4D^{n-\dim(H)}.
\end{equation*}
\end{enumerate}
\end{thm}
The paper is organized as follows. Based on \cite{Rovelli} and \cite{Serre1}, in Section 2 we introduce the biequivariant compactifications  and very ample line bundles that will be used in the remainder of the paper.  In Section 3 we state the definitions and results about derivations needed in this paper. In Section 4 we define and study the function of ideals that are used in the proofs of our main results; it is important to remark that the point of view that we will have is completely algebraically and it wont depend on the exponential map. The proof of our main results follow roughly the standard techniques of the multiplicity estimates theorems, see for instance \cite{Philippon} and \cite{Nakamaye1}; in Section 5 we shall recall and prove the main tools used in these standard techniques. In the last two sections we demonstrate our main statements.
 
\subsection*{Notation and conventions}
In this paper $\field$ is a an algebraically closed subfield of $\Com$ and all the schemes are defined over $\Spec(\field)$ unless otherwise is mentioned. The ring of polynomials (resp.  field of rational functions) in  variables $x_0,\ldots, x_N$ with coefficients in $\field$  is denoted by $\field[x_0,\ldots,x_N]$ (resp. $\field(x_0,\ldots,x_N)$); an element $\frac{f}{h}\in\field(x_0,\ldots,x_N)$ may be considered as a function outside $\Z(h):=\big\{z\in\field^{N+1}:\;h(z)=0\big\}$ and its value in $z\in \field^{N+1}\setminus\Z(h)$ will be denoted by $\frac{f}{h}(x_0,\ldots,x_N)|_{(x_0,\ldots,x_n)=z}$. In this paper $G=(G,\mu,1)$ is  a $n-$dimensional connected  algebraic group  over $\field$  with $\gfr$ its associated Lie algebra.  If there is not possibility of confusion, for all $g,h\in G$ we denote  its product  by $gh$ or $g\cdot h$ instead of $\mu(g,h)$.  The maximal connected affine subgroup of $G$ will be denoted by $L$  which has dimension $m$. If $V$ is a variety, $\Oc_V$ denotes its structural sheaf and $\Oc_{x,V}$ its stalk in $x\in V$.  We call $\field[V]$  the ring of regular functions and $\field(V)$ its function field. Also $\Pic(V)$ denotes  the group of isomorphism classes of line bundles; thus $\Pic(V)$ may be identified with the isomorphism classes of invertible sheaves, see \cite[p. 128]{Hartshorne}. We identify $V(\Com):=\mathrm{Hom}_{\Spec(\field)}(\Spec(\Com),V)$  with the set of closed points of $V\times_{\Spec(\field)}\Spec(\Com).$ If $\Li\in \Pic(V)$ is  base point free, denote by $\phi_\Li:V\rightarrow \Pro^{\dim(\mathrm{H}^0(V,\Li))-1}$ the unique morphism, up to an automorphism of $\Pro^{\dim(\mathrm{H}^0(V,\Li))-1}$, induced by $\Li$. Let $\mathcal{F}$ and $\mathcal{G}$ be  $\Oc_V$-modules. For any $\Oc_V-$modules morphism $\varphi:\mathcal{F}\rightarrow\mathcal{G}$ and $U\subseteq V$ open subset, $\varphi(U):\mathcal{F}(U)\rightarrow\mathcal{G}(U)$ denotes the $\Oc_V(U)-$module morphism. For  subsets $X$ and $Y$ of $\field[x_0,\ldots,x_N]$, $\;^{h}X$ is the set of homogeneous  elements of $X$ and $(X,Y)$ is the ideal generated by $X$ and $Y$. If $I$ and $J$ are ideals  of $\field[x_0,\ldots,x_N]$, $IJ$ (sometimes written $I\cdot J$) is the ideal generated by the product set.  For any $z\in V$, denote by $T_z(V)$ the tangent space of $V$ in $z$.    We write the superindexes in parenthesis if confusion is possible with exponentiation.

\section{Biequivariant compactifications}
The purpose of this section is to construct a particular biequivariant compactification of $G$ and  a projective embedding of it.   
 
\subsection{Construction of $\G$}
 In this subsection we recall the  compactification of $G$  constructed by Rovelli \cite{Rovelli} and we show that it is $G-$biequivariant. The Chevalley-Rosenlicht-Barsotti Theorem, see \cite[Thm. 1.1.1]{Rovelli}, states that $L$ is a normal subgroup of $G$ and  $A:=G/L$ has an abelian variety structure such that the following sequence is an exact sequence of algebraic groups
\begin{equation*}
0\longrightarrow L\longrightarrow G\stackrel{\pi_G}\longrightarrow A\longrightarrow 0.
\end{equation*}
 We say that an algebraic subgroup $H$ of $G$  is \emph{anti-affine} if $\field[H]=\field$. If $K$ is the  smallest algebraic subgroup  of $G$ such that $G/K$ is an affine group, then $K$ is anti-affine, connected and contained in the center of $G$ by \cite[Thm 1.2.1]{BSU}; in this section $K$ denotes the smallest algebraic subgroup of $G$ such that $G/K$ is an affine group and we call it the anti-affine part of $G$. The restriction  $\mu|_{K\times L}$ is a surjective morphism of algebraic groups with kernel $\big\{\big(g^{-1},g\big)\in K\times L:\;g\in L\cap K\big\}$ by \cite[Thm. 1.2.4]{BSU}. It is well known, see \cite[Sec. 1.10]{Borel}, the existence of an injective morphism of algebraic groups $L\hookrightarrow \Gl_m$; hence we assume that $L$ is an algebraic subgroup of $\Gl_m$ from now on. Identifying $\Gl_m$ with the complement of the zero set of the determinant polynomial in $\Af^{m^2}$,  we assume that $\Gl_m$ is embedded in $\Af^{m^2}$. Rovelli \cite[Ex. 3.4.2]{Rovelli}  shows that  $\Pro^{m^2}$ is a  left $\Gl_m-$equivariant compactification  with the following open embedding 
\begin{equation}
\label{E1}
\Gl_m\longrightarrow \Pro^{m^2},\qquad (x_{i,j})\mapsto [1:x_{1,1}:x_{1,2}:\ldots:x_{m,m}]
\end{equation}  
 and the extension of the left translation 
\begin{equation*}
\psi_{L}:\Gl_m\times\Pro^{m^2}\rightarrow\Pro^{m^2} 
\end{equation*}
\begin{equation*}
 \psi_L\big((x_{i,j}), [p_0:...:p_{m^2}]\big)=\Bigg[ p_0:\sum_{j=1}^{m} x_{1,j}p_{(j-1)m+1}:...:\sum_{j=1}^{m} x_{m,j}p_{(j-1)m+m}\Bigg].
\end{equation*}
Moreover, (\ref{E1}) is  a biequivariant compactification of $\Gl_m$ with the following extension of the right translation  
\begin{equation*}
\psi_{R}:\Pro^{m^2}\times\Gl_m\rightarrow\Pro^{m^2} 
\end{equation*}
\begin{equation*}
 \psi_R\big([p_0:\ldots:p_{m^2}],(x_{i,j})\big)=\Bigg[p_0:\sum_{j=1}^{m}p_{j} x_{j,1}:\ldots:\sum_{j=1}^{m} p_{m^2-m+j}x_{j,m}\Bigg].
\end{equation*}
From now on we consider $\Gl_m$ embedded in $\Pro^{m^2}$. Let $\Lb$ be the closure of $L$ in $\Pro^{m^2}$. Then $\Lb$ is a biequivariant compactification of $L$ with the actions $\psi_L|_{L\times \Lb}$ and $\psi_R|_{\Lb\times L}$. The following morphism is a left action of $L$ in $G$
\begin{equation*}
L\times G\longrightarrow G,\qquad (l,g)\mapsto gl^{-1};
\end{equation*}
thus $L$ acts on $G\times \Lb$. Define
\begin{equation*}
\G:=\big(G\times \Lb\big)/\sim\; \text{where }(g,z)\sim \big(g',z'\big)\text{ if } \big(g',z'\big)=\big(gl^{-1},\psi_L(l,z)\big)\text{ for some }l\in L
\end{equation*}
with the projection
\begin{equation*}
\rho:G\times \Lb\longrightarrow \G,\qquad \rho(g,z)=[g,z].
\end{equation*}
If no confusion is possible, we simply denote by $[g,z]$ the projection of $(g,z)\in G\times \Lb$ in the set $\G$.  In \cite[Cor. 1.3.8]{Rovelli} it is shown that the following morphism is an open embedding
\begin{equation*}
G\longrightarrow\G,\qquad g\mapsto[g,1]
\end{equation*} 
and the following left action of $G$ in $\G$ is an extension of the left translation action
\begin{equation*}
\varrho_L:G\times \G\longrightarrow \G,\qquad \varrho_L(g,[g',z])=[gg',z].
\end{equation*}
Define  the morphism 
\begin{equation*}
 \tau:(G\times \Lb)\times L\longrightarrow G\times\Lb,\qquad\tau((g,z),l)=\big(g,\psi_R(z,l)\big).
\end{equation*} Then 
there is a morphism $\eta:\G\times L\rightarrow \G$  which makes the following diagram commutative
 \begin{equation*}\begin{CD}
  (G\times \Lb)\times L @>\tau>>G\times \Lb\\\
@V\rho\times id_L VV @VV\rho V\\
\G\times L @>\eta>> \G.
\end{CD}\end{equation*}
 Hence  $\eta$ has the following property
\begin{equation}
\label{E2}
\eta\big(\varrho_L(g,z),l\big)= \varrho_L\big(g,\eta(z,l)\big)\qquad \forall\;z\in\G,\; l\in L,\;g\in G.
\end{equation}
Since $K$ is in the center of $G$, we have that for all $g,g'\in G$ and $h\in K\cap L$
\begin{equation*}
\eta\big(\varrho_L(g,g'),h\big)=gg'h=ghg'=\varrho_L\big(gh,g'\big); 
\end{equation*}
 moreover, inasmuch as  $G$ is dense in $\G$, we have that
\begin{equation}
\label{E3}
\eta\big(\varrho_L(g,z),h\big)=\varrho_L\big(gh,z\big) \qquad \forall\; z\in\G,\; g\in G,\;h\in K\cap L.
\end{equation}
We prove that there exists a right action of $G$ in $\G$ which is a morphism and an extension of the right translation.
\begin{prop}
\label{R5}
The following morphism  is an extension of the the right translation 
\begin{equation*}
\varrho_R:\G\times G\longrightarrow \G,\qquad \varrho_R\big(z,kl\big)=\eta(\varrho_L(k,z),l) \quad\forall\; k\in K,\;l\in L,\;z\in\G.
\end{equation*}
  Moreover
 \begin{equation*}
 \varrho_L\big(g,\varrho_R\big(z,g'\big)\big)=\varrho_R\big(\varrho_L(g,z),g'\big)\qquad\forall\;g,g'\in G,\;z\in\G .
 \end{equation*}
\end{prop}
\begin{proof}
Define 
\begin{equation*}
\widehat{\varrho_R}:\G\times K\times L\longrightarrow \G,\qquad \widehat{\varrho_R}(z,k,l)=\eta(\varrho_L(k,z),l).
\end{equation*}
From (\ref{E3})
\begin{equation*}
\widehat{\varrho_R}(z,k,l)=z\qquad\forall\;z\in G,\,(k,l)\in \ker(\mu|_{K\times L})
\end{equation*}
and the density of $G$ in $\G$ yields 
\begin{equation*}
\widehat{\varrho_R}(z,k,l)=z\qquad\forall\;z\in  \G,\,(k,l)\in \ker(\mu|_{K\times L})
\end{equation*}
This implies that $\widehat{\varrho_R}$ factorizes through $\G\times G$ and consequently $\varrho_R$ is a well defined morphism. Since $\varrho_L$ and $\eta$ are actions and $K$ is in the center of $G$, it follows easily from (\ref{E2}) that $\varrho_R$ is a right action. To conclude the first claim, note that 
\begin{equation*}
\varrho_R\big(g,kl\big)=\eta(\varrho_L(k,g),l)=(kg)l=g(kl)\qquad \forall\; k\in K,\; l\in L,\;g\in G.
\end{equation*}
The second claim follows straightforward from (\ref{E2}) and (\ref{E3}).
\end{proof}
\subsection{Linearized line bundles}
In the last part of this section,  we give very ample line bundles of the compactification constructed in Section 2.1. Let $V$ be  a variety with a left action of $L$ which is also a morphism and let $\Li$ be line bundle of $V$. We say that $\Li$ is \emph{$L-$linearized} if $\Li$ has a left $L-$action which is a morphism that lifts the action of $L$ in $V$; we denote by $\Pic^L(V)$ the set of isomorphism classes of $L-$linearized line bundles of $V$ and consider $\Pic^L(V)$ a subset of $\Pic(V)$. Let $\zeta$ be the left action of $L$ in $\Li$ and $\pi:\Li\rightarrow \Lb$ a $L-$linearized line bundle. Define 
 \begin{equation*}
 \overline{\Li}:=(G\times \Li)/\sim\; \text{where }(g,z)\sim \big(g',z'\big)\text{ if } \big(g',z'\big)=\big(gl^{-1},\zeta(l,z)\big)\text{ with }l\in L
 \end{equation*}
 with the projection
 \begin{equation*}
\overline{\pi}: \overline{\Li}\rightarrow \G,\qquad \overline{\pi}([g,z])=[g,\pi(z)].
 \end{equation*}
 It is shown in \cite[Lemma 1.2]{Knop-Lange}\footnote{In \cite{Knop-Lange} it is assumed that $G$ is commutative; nonetheless, the proofs of \cite[Lemma 1.2]{Knop-Lange} and \cite[Lemma 1.4]{Knop-Lange} also work without this assumption.}  that $\overline{\pi}:\overline{\Li}\rightarrow \G$ is indeed a line bundle and this induces a function $\Pic^L(\Lb)\rightarrow \Pic(\G)$; moreover, from \cite[Lemma 1.4]{Knop-Lange}
\begin{equation*}
\overline{\Li\otimes \Li'}=\overline{\Li}\otimes \overline{\Li'} \qquad\forall\;\Li,\Li'\in \Pic^L(\Lb).
\end{equation*}

 It is shown in \cite[Lemma 3.4.3]{Rovelli} that $\mathcal{O}(1)\in \Pic^{\Gl_m}\big(\Pro^{m^2}\big)$; hence for all $k\in\Nat$ we have that  $\mathcal{O}(k)\in \Pic^{\Gl_m}\big(\Pro^{m^2}\big)$ and therefore $\Oc(k)|_\Lb\in \Pic^L(\Lb)$. As a consequence of \cite[Sec. 1.3]{Rovelli}, the following  projection is a morphism
  \begin{equation*}
 \pi_{\G}:\G\rightarrow A,\qquad \pi_{\G}([g,x])=\pi_G(g).
 \end{equation*} 
Denote  by $\Pic_R(\G)$   the subset of $\Pic(\G)$ whose elements are the classes
\begin{equation*}
\Big( \pi_{\G}^*\Li_A\otimes \overline{\Oc(1)|_\Lb}\Big)^{\otimes k} = \pi_{\G}^*\Li^{\otimes k}_A\otimes  \overline{\Oc(k)|_\Lb}
\end{equation*}
where $k\in\Nat$ and  $\Li_A\in\Pic(A)$ is  very ample and symmetric. The first important property of $\Pic_R(\G)$ that we will note is that its elements  are very ample and moreover there exist elements of $\Pic_R(\G)$ such that the image of $\G$ under their corresponding  embeddings are normal projective varieties.
\begin{prop}
\label{R6}For any $\Li\in \Pic_R(\G)$, $\Li$ is very ample. Furthermore,  $\phi_{\Li^{\otimes k}}(\G)$  is a projective normal variety if $k>n$.
\end{prop}
\begin{proof}
See \cite[Thm. 2.3.8 and Cor. 3.2.3]{Rovelli}.
\end{proof}
The following result is obtained with the same the ideas of \cite[Lemma 1]{Masser-Wustholz1}.
\begin{prop}
\label{R7}
Let $\phi:\G\rightarrow \Pro^{N}$ be a closed embedding such that $x_k(\phi(1))\neq 0$ for all $k\in\{0,\ldots,N\}$. There exists $c_5\in \Nat$ with the following properties.
\begin{enumerate}
\item[i)]
  There are an open affine finite covering $\big\{U^{(L)}_\alpha\big\}_{\alpha\in \Ac_L}$ of $G\times \G$ and  bihomogeneous polynomials $T^{(L)}_{0,\alpha},\ldots,T^{(L)}_{N,\alpha}\in\field[x_0\ldots, x_N,y_0,\ldots, y_N]$ of bidegree  $(c_5,c_5)$  for all $\alpha\in \Ac_L$ such that
\begin{equation*}
\phi(\varrho_L(g,z))=\Big[T^{(L)}_{0,\alpha}(\phi(g),\phi(z)):\ldots:T^{(L)}_{N,\alpha}(\phi(g),\phi(z))\Big]\qquad  \forall\; (g,z)\in U^{(L)}_\alpha
\end{equation*}
and 
\begin{equation*}
U^{(L)}_\alpha\cap\{(1,z):\;z\in\G\}\neq \emptyset\qquad  \forall\;\alpha\in\Ac_L.
\end{equation*}
\item[ii)]
  There are an open affine finite covering $\big\{U^{(R)}_\alpha\big\}_{\alpha\in \Ac_R}$ of $\G\times G$ and  bihomogeneous polynomials $T^{(R)}_{0,\alpha},\ldots,T^{(R)}_{N,\alpha}\in\field[x_0\ldots, x_N,y_0,\ldots, y_N]$ of bidegree  $(c_5,c_5)$ for all $\alpha\in \Ac_R$ such that
\begin{equation*}
\phi(\varrho_R(z,g))=\Big[T^{(R)}_{0,\alpha}(\phi(z),\phi(g)):\ldots:T^{(R)}_{N,\alpha}(\phi(z),\phi(g))\Big]\qquad  \forall\; (z,g)\in U^{(R)}_\alpha
\end{equation*}
and 
\begin{equation*}
U^{(R)}_\alpha\cap\{(z,1):\;z\in\G\}\neq \emptyset\qquad  \forall\;\alpha\in\Ac_R.
\end{equation*}
\end{enumerate}
\end{prop}
\begin{proof}
 Call $\Af_k:=\{[x_0:\ldots:x_N]\in\Pro^N:\;x_k\neq 0\}$ for $k\in\{0,\ldots,N\}$. Since $G\times \G$ and $\G\times G$ are quasicompact topological spaces, there exist open affine finite coverings $\big\{W^{(L)}_\alpha\big\}_{\alpha\in \Ac_L}$, $\big\{W^{(R)}_\alpha\big\}_{\alpha\in \Ac_R}$ of $G$ and  $\big\{V^{(L)}_\alpha\big\}_{\alpha\in \Ac_L},\big\{V^{(R)}_\alpha\big\}_{\alpha\in \Ac_R}$ of $\G$ with the following properties:
\begin{enumerate}
\item[1)]$1\in \bigcap_{\alpha\in \Ac_L}W^{(L)}_\alpha$ and $1\in \bigcap_{\alpha\in \Ac_R}W^{(R)}_\alpha$.

\item[2)] For all $\alpha,\alpha'\in \Ac_L$ and $\beta,\beta'\in \Ac_R$, there are $ k(\alpha),k(\alpha'),k(\beta),k(\beta')\in\\\{0,\ldots, N\}$ such that $\phi\big(W^{(L)}_\alpha\big)\subseteq \Af_{k(\alpha)}$,  $\phi\big(V^{(L)}_{\alpha'}\big)\subseteq \Af_{k(\alpha')}$, $\phi\big(W^{(R)}_\beta\big)\subseteq \Af_{k(\beta)}$ and  $\phi\big(V^{(R)}_{\beta'}\big)\subseteq \Af_{k(\beta')}$.

\item[3)] Call $U^{(L)}_\alpha:=W^{(L)}_\alpha\times_{\Spec(\field)} V^{(L)}_\alpha$ and $U^{(R)}_\beta:=V^{(R)}_\beta\times_{\Spec(\field)} W^{(R)}_\beta$ for all $\alpha\in \Ac_L, \beta\in \Ac_R$. Then  $\big\{U^{(L)}_\alpha\big\}_{\alpha\in \Ac_L}$  and $\big\{U^{(R)}_\beta\big\}_{\beta\in \Ac_R}$ are open affine finite coverings of $G\times \G$ and $\G\times G$ respectively.

\item[4)] For all $\alpha\in \Ac_L$ and $\beta\in \Ac_R$, there are  $\overline{T^{(L)}_{0,\alpha}},\ldots,\overline{T^{(L)}_{N,\alpha}}\in\field[x_0\ldots, x_N,y_0,\ldots, y_N]$ of bidegree $(c_{6,\alpha},c_{7,\alpha})$ and $\overline{T^{(R)}_{0,\beta}},\ldots,\overline{T^{(R)}_{N,\beta}}\in\field[x_0\ldots, x_N,y_0,\ldots, y_N]$ of bidegree $(c_{8,\beta},c_{9,\beta})$ such that
\begin{align*}
\phi(\varrho_L(g,z))&=\Big[\overline{T^{(L)}_{0,\alpha}}(\phi(g),\phi(z)):\ldots:\overline{T^{(L)}_{N,\alpha}}(\phi(g),\phi(z))\Big]&\forall\; (g,z)\in U^{(L)}_\alpha\\
\phi(\varrho_R(z,g))&=\Big[\overline{T^{(R)}_{0,\beta}}(\phi(z),\phi(g)):\ldots:\overline{T^{(R)}_{N,\beta}}(\phi(z),\phi(g))\Big]&\forall\; (z,g)\in U^{(R)}_\beta.
\end{align*}
\end{enumerate}
Set $c_5:=\max\{c_{6,\alpha},c_{7,\alpha}, c_{8,\beta},c_{9,\beta}:\;\alpha\in \Ac_L,\beta\in \Ac_R\}$. The properties 2) and 3) imply that for each $\alpha\in \Ac_L$ and $\beta\in \Ac_R$ there exist linear forms $F_{W,\alpha}, H_{V,\beta}\in\field[x_0,\ldots,x_N]$ and $F_{V,\alpha}, H_{W,\beta}\in\field[y_0,\ldots,y_N]$ such that 
\begin{align*}
F_{W,\alpha}(\phi(g))\cdot F_{V,\alpha}(\phi(z))\neq 0&\qquad \forall\;(g,z)\in U_\alpha^{(L)}\\
H_{V,\beta}(\phi(z))\cdot H_{W,\beta}(\phi(g))\neq 0&\qquad \forall\;(z,g)\in U_\beta^{(R)}.
\end{align*}
Define the polynomials
\begin{equation*}
T^{(L)}_{k,\alpha}:=\overline{T^{(L)}_{k,\alpha}}\cdot F_{W,\alpha}^{c_5-c_{6,\alpha}}\cdot F_{V,\alpha}^{c_5-c_{7,\alpha}}\qquad\text{ and }\qquad T^{(R)}_{k,\beta}:=\overline{T^{(R)}_{k,\beta}}\cdot H_{V,\beta}^{c_5-c_{8,\beta}}\cdot H_{W,\beta}^{c_5-c_{9,\beta}}
\end{equation*}
for all $k\in\{0,\ldots,N\}$, and note that  i) and ii) are satisfied.
\end{proof}
\begin{rem}
\label{R8} With the notation as in Proposition \ref{R7}, as it is already noted in \cite[Lemma 1]{Masser-Wustholz1}, if $(g,z)\in G\times \G$ is such that not all the values $T^{(L)}_{k,\alpha}(\phi(g),\phi(z))$ are zero for $k\in\{0,\ldots,N\}$, then 
\begin{equation*}
\phi(\varrho_L(g,z))=\Big[T^{(L)}_{0,\alpha}(\phi(g),\phi(z)):\ldots:T^{(L)}_{N,\alpha}(\phi(g),\phi(z))\Big].
\end{equation*}
The equivalent result holds for $\varrho_R$.
\end{rem}

\section{Derivations}
In this section we state the  definition and facts about derivations that will be needed in the next sections. From now on $\G$ is the biequivariant compactification constructed in the previous section, and for simplicity write $gz:=\varrho_{L}(g,z)$ and $zg:=\varrho_{R}(z,g)$ for all $g\in G$ and $z\in\G$. Also we fix $\Li\in\Pic_R(\G)$  and a projective embedding $\phi:=\phi_{\Li^{\otimes n+1}}:\G\rightarrow \Pro^N$; in particular, $\phi(\G)$ is a normal projective variety. We assume without loss of generality that $x_k(\phi(1))\neq 0$ for all $k\in\{0,\ldots,N\}$ and then  $G_k:=\{z\in\G:\;x_k(\phi(z))\neq 0\}$ is a nonempty open subset of $\G$. Call $I(\G)$ the homogeneous prime ideal of $\field[x_0,\ldots, x_N]$ corresponding to $\phi(\G)\subseteq \Pro^N$.  For all $k\in\{0,\ldots,N\}$ we  denote by $I_k(\G)$ the deshomogenization of $I(\G)$ by $x_k$,  and $\phi_k:=\big(\frac{x_0}{x_k},\ldots,\frac{x_N}{x_k}\big)\circ\phi$.  In particular 
\begin{equation*}
\Oc_{\G}(G_k)\cong \field\bigg[\frac{x_0}{x_k},\ldots, \frac{x_N}{x_k}\bigg]\bigg/I_k(\G).
\end{equation*} We denote by $k_z$ the minimum integer $k$ such that $\phi(z)\in G_k$.  When $I$  is a homogeneous ideal of $\field[x_0,\ldots, x_N]$, $\Z(I)$ denotes the zero set of $I$ in $\Pro^N$. 

  For all $g\in G$, let $\lambda_g:G\rightarrow G$ and $\xi_g:\G\rightarrow \G$ be the left translations by $g$ and $\eta_g:\G\rightarrow \G$ be the right translation by $g$. Denote by $\Der(\Oc_G,\Oc_G)$ the set of $\Spec(\field)-$derivations $\Delta:\Oc_G\rightarrow \Oc_G$, and  analogously define $\Der(\Oc_\G,\Oc_\G)$. Recall that the Lie algebra associated to $G$  is  
\begin{equation*}
\gfr:=\big\{\Delta\in \mathrm{Der}(\Oc_G,\Oc_G):\;\Delta\circ\lambda^*_g=\lambda^*_g\circ \Delta\quad\forall\;g\in G\big\}
\end{equation*}
 and set
\begin{equation*} 
  \overline{\gfr}:=\big\{\Delta\in \mathrm{Der}(\Oc_\G,\Oc_\G):\;\Delta\circ\xi^*_g=\xi^*_g\circ \Delta\quad\forall\;g\in G\big\}.
\end{equation*}  
   Name $\mathrm{Der}_\field(\Oc_{1,G},\field)$ the set of $\field-$derivations $\Delta: \Oc_{1,G} \rightarrow \field$  and analogously define $\mathrm{Der}_\field(\Oc_{1,\G},\field)$. Then, see \cite[Expos\'{e} II]{SGA3}, there exist isomorphisms of $\field-$linear spaces
\begin{equation*}
\gfr\cong \mathrm{Der}_\field(\Oc_{1,G},\field)\cong T_1(G)\cong  T_1(\G)\cong \mathrm{Der}_\field(\Oc_{1,\G},\field)\cong \overline{\gfr}.
\end{equation*}
 In particular,   the   restriction map $\overline{\gfr}\rightarrow \gfr$ induced by the open embedding $G\rightarrow \G$  is an isomorphism of $\field-$linear spaces, and from now on we identify $\gfr$ with $\overline{\gfr}$ via this map.
  Fix a basis $\Delta_1,\ldots, \Delta_d$ of $\bfr$ and remember that $\U(\bfr)$ denotes the universal enveloping algebra of $\bfr$. The subset of $\U(\bfr)$ where the elements are of the form  $\Delta_1^{t_1}\ldots \Delta_d^{t_d}$ with $ t_1,\ldots, t_d\in\Nat\cup\{0\}$ is a basis of $\U(\bfr)$, see \cite[Thm. 7.1.9]{Hilgert-Neeb}. Denote by $\U(\bfr,T)$ the $\field$-linear subspace of $\U(\bfr)$ generated by the elements $\Delta_1^{t_1}\ldots \Delta_d^{t_d}$ with $ t_1,\ldots, t_d\in\Nat\cup\{0\}$ and $\sum_{i=1}^dt_i\leq T$. The adjoint representation  is defined as follows
   \begin{equation*}
   \Ad:G\longrightarrow \Gl(\gfr),\qquad \Ad(g)(\Delta)= \eta_g^*\circ\Delta\circ \eta_{g^{-1}}^*.
   \end{equation*}
 If $\Ad(g)(\bfr)\subseteq \bfr$, then the definition of $\Ad(g)(\Delta)$ may be extended for $\Delta\in\U(\bfr)$ in the natural way: $\Ad(g)(\Delta)=\eta_g^*\circ\Delta\circ \eta_{g^{-1}}^*$; in particular $\Ad(g)(\U(\bfr,T))\subseteq \U(\bfr,T)$.
\section{Ideals}
In this section we define the ideals that will be used in the proofs of the main theorems. This is done in the spirit of \cite[Sec. 4]{Philippon}; nonetheless, instead of using the properties of the $d-$parameters subgroups as it is used in the commutative case,  we use that the elements of $\gfr$ are the left invariant derivations. 

Let $I$ be an homogeneous ideal of $\field[x_0,\ldots, x_N]$. We denote by $\mathcal{K}_I$  the set of all homogeneous primary ideals $J$ containing $I$ with $\Z(J)\neq\emptyset$ and 
\begin{equation*}
\Ib (I):= \left\{ \begin{array}{ll}
\field[x_0,\ldots, x_N] & \mbox{if }\mathcal{K}_I=\emptyset\\
&\\
\bigcap_{J\in\mathcal{K}_I}J& \mbox{otherwise}.\end{array} \right.
\end{equation*}
\begin{lem}
\label{R9}
Let $I$ and $J$ be  homogeneous ideal of $\field[x_0,\ldots, x_N]$. If $\Z(J)=\emptyset$, then
\begin{equation*}
\mathcal{K}_I=\mathcal{K}_{I\cap J}=\mathcal{K}_{IJ}.
\end{equation*}
\end{lem}
\begin{proof}
Insomuch as $I\supseteq I\cap J\supseteq IJ$, we get that $\mathcal{K}_I\subseteq\mathcal{K}_{I\cap J}\subseteq\mathcal{K}_{IJ}$ and therefore it is enough to show that $\mathcal{K}_I\supseteq \mathcal{K}_{IJ}$. If  $I\cap \field\neq\{0\}$, then $\mathcal{K}_I=\field[x_0,\ldots, x_N]$ and the statement is true. Thus  we assume that $I\cap \field=\{0\}$. Let $K\in\mathcal{K}_{IJ}$ and $P\in \Ih\setminus\{0\}$. Hence 
\begin{equation*}
(P)J\subseteq IJ\subseteq K
\end{equation*}
and Hilbert's Nullstellensatz, see \cite[Ch. 1]{Hartshorne}, implies that $P\in K$ inasmuch as $\Z(J)=\emptyset$ (in other words, if $P\not\in K$, then $\sqrt{J}\subseteq K$ since $K$ is primary and therefore $\Z(K)\subseteq \Z(J)$ which is impossible). This means that $I\subseteq K$ and thereby $K\in\mathcal{K}_{I}$.
\end{proof}
A straight consequence of the previous lemma is the following statement. 
\begin{cor}
\label{R10}
Let $I$ be a homogeneous ideals of $\field[x_0,\ldots,x_N]$. If $I_1,\ldots,I_r,I_{r+1},\\\ldots,I_k$ are primary homogeneous ideals such that $I=\bigcap_{s=1}^kI_s$ with $\Z(I_1),\ldots,\Z(I_r)$ nonempty and 
\begin{equation*}
\Z(I_{r+1})=\ldots=\Z(I_k)=\emptyset,
\end{equation*}
then $\Ib(I)=\Ib\big(\bigcap_{s=1}^rI_s\big)=\bigcap_{s=1}^rI_s$. In particular, for any  homogeneous ideal $J$ of $\field[x_0,\ldots,x_N]$:
\begin{enumerate}
\item[i)]$\Ib(I\cap J)=\Ib(I)\cap\Ib(J)$.
\item[ii)]$\Ib(I,J)=\big(\Ib(I),\Ib(J)\big)$.
\item[iii)]$\Ib(\Ib(I))=\Ib(I)$.
\item[iv)]$\Z(\Ib(I))=\Z(I)$.
\end{enumerate}
\end{cor}
 The following application of  Hilbert's Nullstellensatz will be used several times.
\begin{prop}
\label{R11}
Let $\{J_\alpha\}_{\alpha\in \Ac}$ and $\{F_\alpha\}_{\alpha\in \Ac}$ be  families of homogeneous ideals of $\field[x_0,\ldots,x_N]$ such that  $\Z(F_\alpha)=\emptyset$ for all $\alpha\in \Ac$. For any  homogeneous ideal $I$ of $\field[x_0,\ldots,x_N]$
\begin{equation*}
\Ib\Bigg(I,\bigcup_{\alpha\in \Ac}J_\alpha\Bigg)=\Ib\Bigg(I,\bigcup_{\alpha\in \Ac}F_\alpha\cdot J_\alpha\Bigg).
\end{equation*}
 \end{prop}
\begin{proof}
The noetherianity of $\field[x_0,\ldots,x_N]$ let us  assume that the index set $\Ac$ is finite. Then  Corollary \ref{R10} ii) yields
\begin{align}
\label{E4}
\Ib\Bigg(I,\bigcup_{\alpha\in \Ac}J_\alpha\Bigg)&=\Bigg(\Ib(I),\bigcup_{\alpha\in \Ac}\Ib(J_\alpha)\Bigg)\nonumber\\
\Ib\Bigg(I,\bigcup_{\alpha\in \Ac}F_\alpha\cdot J_\alpha\Bigg)&=\Bigg(\Ib(I),\bigcup_{\alpha\in \Ac}\Ib(F_\alpha\cdot J_\alpha)\Bigg).
\end{align}
Finally, Lemma \ref{R9} implies that  $\Ib(J_\alpha)=\Ib(F_\alpha\cdot J_\alpha)$ for all $\alpha\in\Ac$, and therefore Corollary \ref{R10} ii) and (\ref{E4}) conclude the proof.
\end{proof} 
 
From Proposition \ref{R7} we fix $c_5\in\Nat$ , the open affine coverings $\big\{U^{(L)}_\alpha\big\}_{\alpha\in \Ac_L}$ of $G\times \G$ and $\big\{U^{(R)}_\alpha\big\}_{\alpha\in \Ac_R}$ of $G\times \G$, and the  bihomogeneous polynomials $T^{(L)}_{0,\alpha},\ldots, T^{(L)}_{N,\alpha},\\ T^{(R)}_{0,\alpha},\ldots,T^{(R)}_{N,\alpha}\in\field[x_0\ldots, x_N,y_0,\ldots, y_N]$ of bidegree $(c_5,c_5)$ such that
\begin{align*}
\phi\big(gz\big)&=\Big[T^{(L)}_{0,\alpha}(\phi(g),\phi(z)):\ldots:T^{(L)}_{N,\alpha}(\phi(g),\phi(z))\Big]\qquad  \forall\; (g,z)\in U^{(L)}_\alpha,\;\alpha\in \Ac_L\nonumber\\
\phi\big(zg\big)&=\Big[T^{(R)}_{0,\alpha}(\phi(z),\phi(g)):\ldots:T^{(R)}_{N,\alpha}(\phi(z),\phi(g))\Big]\qquad  \forall\; (z,g)\in U^{(R)}_\alpha,\;\alpha\in \Ac_R
\end{align*}
and for all $\alpha\in \Ac_L, \beta\in\Ac_R$
\begin{align*}
U^{(L)}_\alpha\cap \{(1,z):\;z\in\G\}\neq \emptyset\quad\text{and}\quad
U^{(R)}_\beta\cap \{(z,1):\;z\in\G\}\neq \emptyset.
\end{align*}
From  Remark \ref{R8} we may assume that for all $\alpha\in \Ac_L$ and $\beta\in \Ac_R$
\begin{align*}
U^{(L)}_\alpha&=G\times\G\setminus\Big\{(g,z)\in G\times\G:\;T^{(L)}_{l,\alpha}(\phi(g),\phi(z))=0\quad\forall\;l\in\{0,\ldots,N\}\Big\}\\
U^{(R)}_\beta&=\G\times G\setminus\Big\{(z,g)\in \G\times G:\;T^{(R)}_{l,\beta}(\phi(z),\phi(g))=0\quad\forall\;l\in\{0,\ldots,N\}\Big\}.
\end{align*} 
 We abbreviate the notation setting   $\x:=(x_0,\ldots,x_N)$, $\y:=(y_0,\ldots,y_N)$ and
\begin{align*}
 T^{(L)}_\alpha(\x,\y):=\Big(T^{(L)}_{0,\alpha}(\x,\y),\ldots,T^{(L)}_{N,\alpha}(\x,\y)\Big)&\qquad \forall\; \alpha\in \Ac_L\nonumber\\
    T^{(R)}_\alpha(\x,\y):=\Big(T^{(R)}_{0,\alpha}(\x,\y),\ldots,T^{(R)}_{N,\alpha}(\x,\y)\Big)&\qquad \forall\; \alpha\in \Ac_R.\nonumber
\end{align*} 
  For  any $f\in \field(x_0,\ldots, x_N)$, set
\begin{align*}
f^{(L)}_{g,\alpha}(\x):=f\Big(T^{(L)}_\alpha(\phi_{k_g}(g),\x)\Big)&\qquad\forall\;\alpha\in \Ac_L\nonumber\\
f^{(R)}_{g,\alpha}(\x):=f\Big(T^{(R)}_{\alpha}(\x,\phi_{k_g}(g))\Big)&\qquad\forall\;\alpha\in \Ac_R
\end{align*}
where recall that $k_g$ is minimum $k$ such that $g\in G_k$. When $I$ is an homogeneous ideal of $\field$, define the ideals
\begin{align*}
\Tb_{L_g}(I)&:=\Ib\Big(P^{(L)}_{g,\alpha},I(\G):\; P\in\Ih, \alpha\in \Ac_L\Big)\\
\Tb_{R_g}(I)&:=\Ib\Big(P^{(R)}_{g,\alpha},I(\G):\; P\in\Ih, \alpha\in \Ac_R\Big)
\end{align*}
\begin{rem}
\label{R12}
Let $I$ be a homogeneous ideal of $\field[x_0,\ldots, x_N]$ and $g\in G$. Corollary \ref{R10} iv) yields 
\begin{equation*}
 \Z(\Tb_{L_g}(I))=\phi\Big(g^{-1}\cdot \phi^{-1}(\Z(I))\Big)\qquad\text{and}\qquad \Z(\Tb_{R_g} (I))=\phi\Big( \phi^{-1}(\Z(I))\cdot g^{-1}\Big).
 \end{equation*}
\end{rem}
We start proving the main properties of the ideals defined above. 
   \begin{lem}
  \label{R13}
 Let $I$ be a homogeneous ideal of $\field[x_0,\ldots, x_N]$ and  $g\in G$. Then
  \begin{equation*}
\Tb_{L_g}(\Ib(I))=\Tb_{L_g}(I)\qquad\text{and}\qquad \Tb_{R_g}(\Ib(I))=\Tb_{R_g}(I).
  \end{equation*}
   \end{lem}
\begin{proof}
We just show the first equality since the second equality is proven analogously. Inasmuch as $I\subseteq \Ib(I)$, we get
\begin{equation*}
\Tb_{L_g}(I)\subseteq \Tb_{L_g}(\Ib(I)).
\end{equation*}
 By Corollary \ref{R10} and the Primary Decomposition Theorem,  there is an homogeneous ideal $J$ of $\field[x_0,\ldots,x_N]$ such that $I=\Ib(I)\cap J$ and $\Z(J)=\emptyset$. Given a homogeneous ideal $K$  of $\field[x_0,\ldots,x_N]$, write
\begin{equation*}
K_g:=\big(P_{g,\alpha}^{(L)},I(\G):\; P\in\,^{h}K,\,\alpha\in\Ac_L\big).
\end{equation*}
 Remark \ref{R12} asserts  that $\Z(\Tb_{L_g}(J))=\emptyset$ and consequently $\Z(J_g)=\emptyset$ by Corollary \ref{R10} iv). Then we complete the proof as follows
\begin{align*}
\Tb_{L_g}(\Ib(I))&=\Ib(\Ib(I)_g)\\
&=\Ib\big(\Ib(I)_g\cdot J_g\big)&\text{by Proposition \ref{R11}}\\
&\subseteq \Ib\big((\Ib(I)\cdot J)_g\big)&\text{since }\Ib(I)_g\cdot J_g=(\Ib(I)\cdot J)_g\\
&\subseteq \Ib(I_g)&\text{since }\Ib(I)\cdot J\subseteq I\\
&=\Tb_{L_g}(I).
\end{align*}
\end{proof}
\begin{lem}
\label{R14}
Let $I$ be  a homogeneous ideal of $\field[x_0,\ldots,x_N]$ and $g,h\in G$.
\begin{enumerate}
\item[i)]$\Tb_{L_h}\big(\Tb_{L_g}(I)\big)=\Tb_{L_{gh}} (I)$.
\item[ii)]$\Tb_{R_h}\big(\Tb_{R_g}(I)\big)=\Tb_{R_{hg}} (I)$.
\item[iii)]$\Tb_{R_h}\big(\Tb_{L_g}(I)\big)=\Tb_{L_g}\big(\Tb_{R_h}(I)\big)$.
\item[iv)]$\Tb_{L_1}(I)=\Ib\big(P,I(\G):\;P \in \Ih\big)=\Tb_{R_1}(I)$.
\end{enumerate}
\end{lem}
\begin{proof}
We only show i) since ii),iii) and iv) are proven in a very similar way. From Lemma \ref{R13} it is enough to show 
\begin{equation}
\label{E5}
\Ib\bigg(\Big(P^{(L)}_{g,\alpha}\Big)^{(L)}_{h,\beta},I(\G):\ P\in \Ih,\, \alpha,\beta\in \Ac_L\bigg)=\Tb_{L_{gh}}(I)
\end{equation}
Given $\alpha,\beta,\gamma\in \Ac_L$, let $U$ be the subset of $\G$ which elements $z$ satisfy  that $(h,z)\in U^{(L)}_\beta$,  $(g,hz)\in U^{(L)}_\alpha$ and $(gh,z)\in U^{(L)}_\gamma$; in particular $U$ is an open nonempty subset. Then for all $z\in U$
\begin{align}
\label{E6}
\phi\big((gh)z\big)&=\Big[T^{(L)}_{0,\gamma}(\phi(gh),\phi(z)):\ldots:T^{(L)}_{N,\gamma}(\phi(gh),\phi(z))\Big]\nonumber\\
&=\Big[T^{(L)}_{0,\alpha}\Big(\phi(g),T^{(L)}_{\beta}(\phi(h),\phi(z))\Big):\ldots:T^{(L)}_{N,\alpha}\Big(\phi(g),T^{(L)}_{\beta}(\phi(h),\phi(z))\Big)\Big]\nonumber\\
&=\phi\big(g(hz)\big).
\end{align}
For all $l\in\{0,\ldots,N\}$ and $P\in\Ih$
\begin{align*}
w^{(l)}_{\alpha,\beta}(\x)&:=T^{(L)}_{l,\alpha}\big(\phi_{k_g}(g),T^{(L)}_{\beta}(\phi_{k_h}(h),
\x)\big)^{\deg(P)}\\
v^{(l)}_\gamma(\x)&:=T^{(L)}_{l,\gamma}(\phi_{k_{gh}}(gh),\x)^{\deg(P)}.
\end{align*}
From (\ref{E6}) we deduce that for all $z\in U$
\begin{align*}
P\Big(T^{(L)}_{\gamma}(\phi_{k_{gh}}(gh),\phi_{k_z}(z))\Big)w^{(l)}_{\alpha,\beta}(\phi_{k_z}(z))&=\\
P\Big(T^{(L)}_{\alpha}\Big(\phi_{k_g}(g),T^{(L)}_{\beta}(\phi_{k_h}(h),\phi_{k_z}(z))\Big)\Big)v^{(l)}_\gamma(\phi_{k_z}(z)).&
\end{align*}
Since $U$ is dense in $\G$
\begin{equation*}
P\Big(T^{(L)}_{\gamma}(\phi_{k_{gh}}(gh),\x)\Big)w^{(l)}_{\alpha,\beta}(\x)-P\Big(T^{(L)}_{\alpha}\Big(\phi_{k_g}(g),T^{(L)}_{\beta}(\phi_{k_h}(h),\x)\Big)\Big)v^{(l)}_\gamma(\x)\in I(\G),
\end{equation*}
 and consequently we have the equality
 \begin{align}
\label{E7}
 \Ib\Big(P\Big(T^{(L)}_{\alpha}\Big(\phi_{k_g}(g),T^{(L)}_{\beta}(\phi_{k_h}(h),\x)\Big)\Big)v^{(l)}_\gamma(\x), I(\G):\, &\nonumber\\
 P\in \Ih,\,\alpha,\beta,\gamma\in\Ac_L,\;  l\in\{0,\ldots,N\}\Big)&=\nonumber\\
 \Ib\Big(P\Big(T^{(L)}_{\gamma}(\phi_{k_{gh}}(gh),\x)\Big)w^{(l)}_{\alpha,\beta}(\x), I(\G):\; P\in \Ih,\; \alpha,\beta,\gamma\in\Ac_L,\; l\in\{0,\ldots,N\}\Big).
\end{align}
On the other hand 
\begin{align*}
\Z\Big(w^{(l)}_{\alpha,\beta}(\x), I(\G):\; l\in\{0,\ldots,N\}, \alpha,\beta\in \Ac_L\Big)=&\\
\Z\Big(v^{(l)}_\gamma(\x),I(\G):\; l\in\{0,\ldots,N\}, \gamma\in \Ac_L\Big)=&\emptyset.
\end{align*}
 Then (\ref{E5}) is a straight consequence of applying  Proposition \ref{R11}  to both sides of (\ref{E7}). 
 \end{proof}

 Recall that $\{\Delta_1,\ldots,\Delta_d\}$ is a fixed basis of $\bfr$. Let $c_6>1$ be a big enough natural number with the following property: for all $l,k\in\{0,\ldots, N\}$ and $j\in\{1,\ldots,d\}$ there is $R^{(l)}_{j,k}\in\field[x_0,\ldots,x_N]$  of total degree $c_6$  such that 
\begin{equation*}
\Delta_j(G_k)\bigg(\frac{x_l}{x_k}+I_k(\G)\bigg)=R^{(l)}_{j,k}\bigg(\frac{x_0}{x_k},\ldots, \frac{x_N}{x_k}\bigg)+I_k(\G)
\end{equation*}
and we define the homogenizations 
\begin{equation*}
Q^{(l)}_{j,k}(x_0,\ldots,x_N):=x_k^{c_6}\cdot R^{(l)}_{j,k}\bigg(\frac{x_0}{x_k},\ldots, \frac{x_N}{x_k}\bigg).
\end{equation*}
For $f\in\field(x_0,\ldots,x_N,y_0,\ldots,y_N)$,  $k\in\{0,\ldots,N\}$ and $j\in\{1,\ldots,d\}$, write
 \begin{align*}
   \Dc(\Delta_j)(f(\x,\y))&:=\sum_{l=0}^N\frac{\partial f(\x,\y)}{\partial y_l}\cdot Q^{(l)}_{j,0}(\y),
 \end{align*}
  and  $\Dc(1)(f)=f$ where $1$ is the multiplicative neutral element of $\U(\bfr)$.  In general, for  $t_1,\ldots ,t_d\in\Nat\cup\{0\}$, set
\begin{align*}
\Dc\big(\Delta^{t_1}_1\ldots \Delta^{t_d}_d\big)(f)&:=\overbrace{\Dc(\Delta_{d})\circ\ldots\circ \Dc(\Delta_{d})}^{t_d}\circ\ldots\circ\overbrace{\Dc(\Delta_{1})\circ\ldots\circ \Dc(\Delta_{1})}^{t_1}(f),
\end{align*}
and then the definition of  $\Dc(\Delta)(f)$ for $\Delta\in\U(\bfr)$ is extended  by linearity. Let $P\in\field[x_0,\ldots, x_N]$ be homogeneous, $\alpha\in\Ac_R$ and  $\Delta\in\U(\bfr,T)$. Define
\begin{align*}
P_{\Delta,\alpha}(\x)&:=\Dc(\Delta)\Big(P\Big(T^{(R)}_\alpha(\x,\y)\Big)\Big)\Big|_{\y=\phi_0(1)}\in \field[x_0,\ldots, x_N].
\end{align*} 
 Set
\begin{align*}
\partial_{L_g}^T(I)&:=\bigg(\big(P_{\Delta,\alpha}\big)^{(L)}_{g,\beta},I(\G):\; P\in \Ih,\;\alpha\in\Ac_R,\;\beta\in \Ac_L,\;\Delta\in \U(\bfr,T)\bigg)\nonumber\\
 \partial_{R_g}^T(I)&:=\bigg(\big(P_{\Delta,\alpha}\big)^{(R)}_{g,\beta},I(\G):\; P\in \Ih,\;\alpha,\beta\in\Ac_R,\;\Delta\in \U(\bfr,T)\bigg).
\end{align*}
\begin{rem}
\label{R15}
Let $I$ be  homogeneous ideal of $\field[x_0,\ldots,x_N]$, $T\in\Nat\cup \{0\}$ and $g\in G$. If $I(\G)$ and $I$ are generated by homogeneous  polynomials of degree at most $D$, then $\partial_{L_g}^T(I)$ and $\partial_{R_g}^T(I)$ will be generated by homogeneous polynomials of degree at most $c_5^2D$; in particular, the previous upper bound  is independent of $T$.
\end{rem}
The main goal of this section is to show that $\Ib(\partial_{L_g}^T(I))$ and $\Ib(\partial_{R_g}^T(I))$ have similar properties to the ideals defined in \cite[D\'{e}f. 4.2]{Philippon}. Thus the remainder of this section is devoted to this goal. To achieve this  
purpose, we need some technical results and, to state these auxiliary lemmas, we need some definitions. For   $e\in\field(x_0,\ldots,x_N)$, $k\in\{0,\ldots,N\}$ and $j\in\{1,\ldots,d\}$, write
 \begin{align*}
  \Bc_k(\Delta_j)(e(\x))&:=\sum_{l=0}^N\frac{\partial e(\x)}{\partial x_l}\cdot Q^{(l)}_{j,k}(\x)
  \end{align*}
  and  $\Bc_k(1)(e)=e$  where $1$ is the multiplicative neutral element of $\U(\bfr)$.  We extend the definition of  $\Bc_k(\Delta)$ for  $\Delta\in\U(\bfr)$ as follows: for  $t_1,\ldots ,t_d\in\Nat\cup\{0\}$ set
\begin{align*}
\Bc_k\big(\Delta^{t_1}_1\ldots \Delta^{t_d}_d\big)(e)&:=\overbrace{\Bc_k(\Delta_{d})\circ\ldots\circ \Bc_k(\Delta_{d})}^{t_d}\circ\ldots\circ\overbrace{\Bc_k(\Delta_{1})\circ\ldots\circ \Bc_k(\Delta_{1})}^{t_1}(e)
\end{align*}
and then $\Bc_k(\Delta)(e)$ is  extended  by linearity for $\Delta\in\U(\bfr)$. Let $P\in\field[x_0,\ldots, x_N]$ be homogeneous, $\alpha\in\Ac_R$, $\Delta\in\U(\bfr,T)\setminus\U(\bfr,T-1)$, $\Delta_1',\ldots,\Delta'_r\in\bfr$ and $ k,k_1,\ldots, k_r\in\{0,\ldots, N\}$. Define the polynomials
\begin{align*}
P^{(\Delta'_1,\ldots,\Delta'_r)}_{(k_1,\ldots,k_r)}(\x)&:=\Bc_{k_r}(\Delta'_r)\circ\ldots\circ\Bc_{k_1}(\Delta'_1)(P(\x))\\
P^{\Delta}_{k}(\x)&:=\Bc_k(\Delta)(P(\x))\\
P^{(k)}_{\Delta,\alpha}(\x)&:=T^{(R)}_{k,\alpha}(\x,\y)^{\deg(P)+T}\cdot\Dc(\Delta)\Bigg(\frac{P\big(T^{(R)}_\alpha(\x,\y)\big)}{T^{(R)}_{k,\alpha}(\x,\y)^{\deg(P)}}\Bigg)\Bigg|_{\y=\phi_0(1)}.
\end{align*} 
Since $\Delta_j(G_k)(I_k(\G))\subseteq I_k(\G)$ for all $j\in\{1,\ldots,d\}$ and $k\in\{0,\ldots, N\}$, it should be clear that if $P\in I(\G)$ is homogeneous, then, with the notation as above, $P^{(\Delta'_1,\ldots,\Delta'_r)}_{(k_1,\ldots,k_r)}, P^{\Delta}_{k}, P^{(k)}_{\Delta,\alpha}$ and $ P_{\Delta,\alpha}$ are also in $I(\G)$.  Let $I$ be a homogeneous ideal of $\field[x_0,\ldots, x_N]$ and $T\in\Nat\cup\{0\}$. Call
 \begin{align*}
    \Bb^T(I)&:=\Ib\bigg(P^{(\Delta'_1,\ldots,\Delta'_r)}_{(k_1,\ldots,k_r)},I,I(\G):\; P\in \Ih,\,r\in\{1,\ldots,T\},\\&\qquad\qquad\Delta'_1,\ldots,\Delta'_r\in\bfr,\, k_1,\ldots,k_r\in\{0,\ldots, N\}\bigg)\\
    \Cb^T(I)&:=\Ib\Big(P^{\Delta}_{k},I(\G):\; P\in \Ih,\;\Delta\in\U(\bfr,T),\;k\in\{0,\ldots, N\}\Big)\\
 \Db^T(I)&:=\Ib\Big(P^{(k)}_{\Delta,\alpha}, I(\G):\; P\in \Ih,\; \alpha\in \Ac_R,\;\Delta\in\U(\bfr,T),\;k\in\{0,\ldots, N\}\Big)\\
 \Eb^T(I)&:=\Ib\Big(P_{\Delta,\alpha},I(\G):\; P\in \Ih,\; \alpha\in \Ac_R,\;\Delta\in\U(\bfr,T)\Big).
  \end{align*}
  We shall show that $\Bb^T(I)=\Cb^T(I)=\Db^T(I)=\Eb^T(I)$. We start proving an auxiliary statement.
\begin{lem}\label{R16}
Let $P\in\field[x_0,\ldots, x_N]$ be homogeneous, $\Delta\in\bfr$, $\alpha,\beta\in\Ac_R$ and $ k,l\in\{0,\ldots,N\}$. Set $u(\x):=\frac{(x_k^{\deg(P)})^{(l)}_{\Delta,\alpha}}{T^{(R)}_{k,\alpha}(\x,\phi_0(1))^{\deg(P)-1}}$. Then the following polynomials are in $I(\G)$. 
\begin{enumerate}
\item[i)]
$x_k^{c_6-1}\cdot P^\Delta_l(\x)-x_l^{c_6-1}\cdot P^\Delta_k(\x)-\deg(P)\cdot \big(x_k\big)_l^\Delta\cdot x_k^{c_6-2}\cdot P(\x)$
\item[ii)]
$T^{(R)}_{k,\alpha}(\x,\phi_0(1))\cdot P^{(l)}_{\Delta,\alpha}(\x)
-T^{(R)}_{l,\alpha}(\x,\phi_0(1))\cdot P^{(k)}_{\Delta,\alpha}(\x)
- u(\x)\cdot P\big(T^{(R)}_\alpha(\x,\phi_0(1))\big)$.
\item[iii)]$T^{(R)}_{k,\alpha}(\x,\phi_0(1))^{\deg(P)}\cdot P_{\Delta,\beta}(\x)-T^{(R)}_{k,\beta}(\x,\phi_0(1))^{\deg(P)}\cdot P_{\Delta,\alpha}(\x)\\+\big(x_k^{\deg(P)}\big)_{\Delta,\alpha}\cdot P\big(T^{(R)}_{\beta}(\x,\phi_0(1))\big)-\big(x_k^{\deg(P)}\big)_{\Delta,\beta}\cdot P\big(T^{(R)}_{\alpha}(\x,\phi_0(1))\big)$
\end{enumerate}
\end{lem}
\begin{proof}
 First we show  that i) is in $I(\G)$. For all $z\in G_k\cap G_l$, the Leibniz's rule yields
\begin{align}
\label{E8}
\Delta(G_l)\bigg(\frac{P(\x)}{x_l^{\deg(P)}}\bigg)\bigg|_{\x=\phi_l(z)}&=\Delta(G_l\cap G_k)\bigg(\frac{P(\x)}{x_l^{\deg(P)}}\bigg)\bigg|_{\x=\phi_l(z)}\nonumber\\
&=\Delta(G_l\cap G_k)\bigg(\frac{P(\x)}{x_k^{\deg(P)}}\cdot\frac{x_k^{\deg(P)}}{x_l^{\deg(P)}}\bigg)\bigg|_{\x=\phi_l(z)}\nonumber\\
&=\Delta(G_k)\bigg(\frac{P(\x)}{x_k^{\deg(P)}}\bigg)\cdot\frac{x_k^{\deg(P)}}{x_l^{\deg(P)}}\nonumber\\
&\quad+\Delta(G_l)\bigg(\frac{x_k^{\deg(P)}}{x_l^{\deg(P)}}\bigg)\cdot\frac{P(\x)}{x_k^{\deg(P)}}\bigg|_{\x=\phi_l(z)}.
\end{align}
Multiplying (\ref{E8}) by $x_l^{\deg(P)+c_6-1}\cdot x_k^{c_6-1}$, we get 
\begin{align}
\label{E9}
x_k^{c_6-1}\cdot P^\Delta_l(\x)\Big|_{\x=\phi_l(z)}&=
x_l^{c_6-1}\cdot P^\Delta_k(\x)+\big(x_k^{\deg(P)}\big)_l^\Delta\cdot x_k^{c_6-1}\cdot \frac{P(\x)}{x_k^{\deg(P)}}\Big|_{\x=\phi_l(z)}\nonumber\\
&=x_l^{c_6-1}\cdot P^\Delta_k(\x)+\deg(P)\cdot \big(x_k\big)_l^\Delta\cdot x_k^{c_6-2}\cdot P(\x)\Big|_{\x=\phi_l(z)}.
\end{align}
 Insomuch as $c_6>1$, $x_k^{c_6-2}$ is certainly a polynomial.  Since $G_k\cap G_l$ is dense in $\G$, (\ref{E9}) implies that i) evaluated in any point of $\phi(\G)$ is zero and therefore it is contained in $I(\G)$.

Now we show that ii) is contained in $I(\G)$. For $0\leq j\leq N$ and $w\in G_0$, set 
\begin{align*}
P^{(j,w)}_{\Delta,\alpha}(\x)&:=T^{(R)}_{j,\alpha}(\x,\y)^{\deg(P)+1}\cdot\Dc(\Delta)\Bigg(\frac{P\big(T^{(R)}_\alpha(\x,\y)\big)}{T^{(R)}_{j,\alpha}(\x,\y)^{\deg(P)}}\Bigg)\Bigg|_{\y=\phi_0(w)}; 
\end{align*}
hence $P^{(j,1)}_{\Delta,\alpha}=P^{(j)}_{\Delta,\alpha}$. Let $(z,w)$ be an element of $U^{(R)}_\alpha$ such that $z\in G_l\cap G_k$,  $w\in G_0$ and $zw\in G_l\cap G_k$. Call $v:=(\phi_l(z),\phi_0(w))$. Then
\begin{equation*}
\frac{P\big(T^{(R)}_\alpha(\x,\y)\big)}{T^{(R)}_{l,\alpha}(\x,\y)^{\deg(P)}} \Bigg|_{(\x,\y)=v}=\frac{P\big(T^{(R)}_\alpha(\x,\y)\big)}{T^{(R)}_{k,\alpha}(\x,\y)^{\deg(P)}}\cdot \frac{T^{(R)}_{k,\alpha}(\x,\y)^{\deg(P)}}{T^{(R)}_{l,\alpha}(\x,\y)^{\deg(P)}} \Bigg|_{(\x,\y)=v}
\end{equation*}
and  Leibniz's rule yields
\begin{align}
\label{E10}
\Dc(\Delta)\bigg(\frac{P\big(T^{(R)}_\alpha(\x,\y)\big)}{T^{(R)}_{l,\alpha}(\x,\y)^{\deg(P)}}\bigg) \Bigg|_{(\x,\y)=v}
&=\Dc(\Delta)\bigg(\frac{P\big(T^{(R)}_\alpha(\x,\y)\big)}{T^{(R)}_{k,\alpha}(\x,\y)^{\deg(P)}}\bigg)\cdot \frac{T^{(R)}_{k,\alpha}(\x,\y)^{\deg(P)}}{T^{(R)}_{l,\alpha}(\x,\y)^{\deg(P)}}\nonumber\\
&+\Dc(\Delta)\bigg(\frac{T^{(R)}_{k,\alpha}(\x,\y)^{\deg(P)}}{T^{(R)}_{l,\alpha}(\x,\y)^{\deg(P)}}\bigg)\cdot \frac{P\big(T^{(R)}_\alpha(\x,\y)\big)}{T^{(R)}_{k,\alpha}(\x,\y)^{\deg(P)}}\bigg|_{(\x,\y)=v}
\end{align}
Multiplying (\ref{E10}) by $T^{(R)}_{l,\alpha}(\x,\y)^{\deg(P)+1}\cdot T^{(R)}_{k,\alpha}(\x,\y)$, we arrive to
\begin{align}
\label{E11}
T^{(R)}_{k,\alpha}(\x,\y)\cdot P^{(l,\y)}_{\Delta,\alpha}(\x)\Big|_{(\x,\y)=v}
&=T^{(R)}_{l,\alpha}(\x,\y)\cdot P^{(k,\y)}_{\Delta,\alpha}(\x)\nonumber\\
&\qquad+T^{(R)}_{k,\alpha}(\x,\y)\cdot (x_k^{\deg(P)})^{(l,\y)}_{\Delta,\alpha}\cdot \frac{P\big(T^{(R)}_\alpha(\x,\y)\big)}{T^{(R)}_{k,\alpha}(\x,\y)^{\deg(P)}}\bigg|_{(\x,\y)=v}.
\end{align}
Now note that the set of elements $(z,w)$ in $U^{(R)}_\alpha$ such that $z\in G_l\cap G_k$,  $w\in G_0$ and $zw\in G_l\cap G_k$ is dense in $\G\times G$ since it is intersection of the nonempty open sets
\begin{equation*}
U^{(R)}_\alpha\cap (\G\times G_0)\cap \big((G_l\cap G_k)\times G\big)\cap\{(z,w)\in \G\times G:\;zw\in G_l\cap G_k\}.
\end{equation*}
Thus (\ref{E11}) holds true for all $(z,w)\in \G\times G$ and it leads to ii) taking $w=1$.

We prove iii). Let $(z,w)$  be in $U^{(R)}_\alpha\cap U^{(R)}_\beta$ and write $v:=(\phi_{k_z}(z), \phi_{k_w}(w))$. Note that
\begin{align}
\label{E12}
T^{(R)}_{k,\alpha}(\x,\y)^{\deg(P)}\cdot P\Big(T^{(R)}_{\beta}(\x,\y)\Big)\bigg|_{(\x,\y)=v}&=\nonumber\\
T^{(R)}_{k,\beta}(\x,\y)^{\deg(P)}\cdot P\Big(T^{(R)}_{\alpha}(\x,\y)\Big)\bigg|_{(\x,\y)=v}
\end{align}
since
\begin{align*}
\phi(zw)=\Big[T^{(R)}_{0,\alpha}(z,w):\ldots:T^{(R)}_{N,\alpha}(z,w)\Big]=\Big[T^{(R)}_{0,\beta}(z,w):\ldots:T^{(R)}_{N,\beta}(z,w)\Big];
\end{align*}
moreover, (\ref{E12}) remains valid also for $(z,w)\not\in U^{(R)}_\alpha\cap U^{(R)}_\beta$ since $T^{(R)}_{l,\alpha}(z,w)=0$ for all $l\in\{0,\ldots,N\}$ or $T^{(R)}_{l,\beta}(z,w)=0$ for all $l\in\{0,\ldots,N\}$. Applying $\Dc(\Delta)$ to (\ref{E12}), the Leibniz's rule gives  iii) taking $w=1$.
\end{proof}
\begin{lem}
 \label{R17}
Let $I$ be a homogeneous ideal of $\field[x_0,\ldots,x_N]$ and  $T\in \Nat\cup\{0\}$. Then
\begin{equation*}
\Bb^{T}(I)=\Cb^{T}(I).
\end{equation*}
 \end{lem}
 \begin{proof}
The proof is by induction on $T$.  The result is trivial for $T\in\{0,1\}$, and hence we assume that $\Bb^{T-1}(I)=\Cb^{T-1}(I)$ and $T\geq 2$. In the induction step, the inclusion $\Bb^T(I)\supseteq \Cb^T(I)$ is trivial since 
\begin{equation*}
P^{(\Delta'_1,\ldots,\Delta'_r)}_{(k,\ldots,k)}=P^{\Delta'_1\ldots\Delta'_r}_k\qquad\forall\; P\in\Ih,\;r\in\{0,\ldots,T\},\;k\in\{0,\ldots,N\},\;\Delta'_1\ldots\Delta'_r\in\bfr.
\end{equation*}
Hence, to complete the proof, it is enough to show that $\Bb^T(I)\subseteq \Cb^T(I)$. For  $s,t,k\in\Zet$ with $k\in\{0,\ldots, N\}$ and $s\in\{0,\ldots,T\}$, define the index sets
\begin{align*}
\Ac_{t,s}&:=\Big\{(k_1,..., k_r)\in\Zet^r:\;r\in\{1,...,t\},\;k_1,...,k_r\in\{0,..., N\},\; k_1=...=k_{\min\{r,s\}}\Big\},
\end{align*}
the monomials  $w_{s,k}:=x^{(T-s)(c_6-1)}_{k}$ and the ideals
\begin{align*}
I_{s}&:=\Ib\bigg( w_{s,k'_0}\cdot P^{(\Delta'_1,\ldots,\Delta'_r)}_{(k'_1,\ldots,k'_r)}, P^{(\Delta'_1,\ldots,\Delta'_r)}_{(k_1,\ldots,k_r)}, I, I(\G):\; P\in \Ih,\;r\in\{1,\ldots,T\},\\
&\qquad\qquad k'_0=k'_{\min\{r,s\}},\; \Delta'_1,...,\Delta'_r\in\bfr,\;(k'_1,...,k'_r)\in \Ac_{T,s},\; (k_1,...,k_r)\in \Ac_{T-1,s}\bigg).
\end{align*}
Take $P\in \Ih$, $k\in\{0,\ldots,N\}$, $s,r\in\{0,\ldots,T\}$ with $s<r$, $\Delta'_1,\ldots,\Delta'_r\in\bfr$ and $(k'_1,\ldots,k'_r)\in \Ac_{T,s}$. From Lemma \ref{R16} i), we get the inclusion
\begin{align}
\label{E13}
\bigg(x_{k'_0}^{c_6-1}\cdot P^{(\Delta'_1,\ldots,\Delta'_s,\Delta'_{s+1})}_{(k'_1,\ldots,k'_s,k_{s+1}')},I(\G)\bigg)\subseteq& \bigg( P^{(\Delta'_1,\ldots,\Delta'_s,\Delta'_{s+1})}_{(k'_1,\ldots,k'_s,k_0')},P^{(\Delta'_1,\ldots,\Delta'_s)}_{(k'_1,\ldots,k'_s)},I(\G)\bigg).
\end{align} Thus from (\ref{E13}) we get that 
\begin{equation*}
\bigg(w_{s+1,k'_0}\cdot P^{(\Delta'_1,\ldots,\Delta'_{s+1})}_{(k'_1,\ldots,k'_{s+1})}\bigg)^{(\Delta'_{s+2},\ldots,\Delta'_{r})}_{(k'_{s+2},\ldots,k'_{r})}\in I_{s+1};
\end{equation*}
 hence by  the Leibniz's rule 
 \begin{align*}
w_{s+1,k'_0}\cdot P^{(\Delta'_1,\ldots,\Delta'_r)}_{(k'_1,\ldots,k'_r)}\in I_{s+1}.
\end{align*}
and consequently
\begin{align}
\label{E14}
w_{s,k'_0}\cdot P^{(\Delta'_1,\ldots,\Delta'_r)}_{(k'_1,\ldots,k'_r)}\in I_{s+1}.
\end{align}
On the other hand, define the following ideals for $s\in\{1,\ldots, T-1\}$
\begin{align*}
J_s&:=\Ib\bigg( P^{(\Delta'_1,\ldots,\Delta'_r)}_{(k_1,\ldots,k_r)}, I, I(\G):\; P\in \Ih,\;r\in\{1,\ldots T-1\},\\
&\qquad\qquad  \Delta'_1,\ldots,\Delta'_r\in\bfr,\;(k_1,\ldots,k_r)\in \Ac_{T-1,s}\bigg)\\
&\subseteq I_s\end{align*} and note that
\begin{equation*}
\Bb^{T-1}(I)=J_1\supseteq\ldots\supseteq J_{T-1}=\Cb^{T-1}(I).
\end{equation*}
Thus the hypothesis of induction establishes that $J_s=J_{s+1}$ for all $s\in\{1,\ldots, T-2\}$; in particular, for all $P\in \Ih$, $s,r\in\{1,\ldots,T-1\}$, $\Delta'_1,\ldots,\Delta'_r\in\bfr$ and $(k_1,\ldots,k_r)\in \Ac_{T-1,s}$
\begin{align}
\label{E15}
 P^{(\Delta'_1,\ldots,\Delta'_r)}_{(k_1,\ldots,k_r)}\in J_s=J_{s+1}\subseteq  I_{s+1}.
\end{align}
Corollary \ref{R10} ii) and the inclusions (\ref{E14}) and (\ref{E15}) let us conclude that 
\begin{align}
\label{E16}
I_0\subseteq I_1\subseteq\ldots\subseteq I_T=\Cb^T(I)
\end{align}
Since $\Z(w_{0,k}:\;k\in\{0,\ldots, N\})=\emptyset$, Proposition \ref{R11} implies that $I_0=\Bb^T(I)$ and then the induction is completed by (\ref{E16}).
 \end{proof}
 For all $T\in\Nat$, define
 \begin{equation*}
 \U_T:=\bigg\{\Delta_1^{t_1}\ldots\Delta_d^{t_d}:\;t_1,\ldots,t_d\in\Nat\cup\{0\},\;\sum_{i=1}^dt_i=T\bigg\}.
 \end{equation*}
 and $\U_0=\{1\}$ (here $1$ is the neutral element of $\U(\bfr)$).
\begin{lem}
\label{R18}
Let $I$ be a homogeneous ideal of $\field[x_0,\ldots,x_N]$ and  $T\in \Nat\cup\{0\}$. Then
\begin{equation*}
\Cb^{T}(I)=\Db^{T}(I).
\end{equation*}
\end{lem}
\begin{proof}
 The proof is by induction on $T$; however, before we start the induction, we show the main tool that relates $\Cb^{T}(I)$ and $\Db^{T}(I)$. Let $g$ be an element of $G_k$ for some $k\in\{0,\ldots,N\}$ and $U:=g^{-1}G_k\cap G_0$. Recall that the elements of $\bfr$ are left invariant derivations; hence for all $ \Delta'_1,\ldots,\Delta'_r\in\bfr$, the following diagram commutes
\begin{equation}
\label{E17}
\begin{CD}
  \Oc_\G(G_k) @>\Delta'_1(G_k)>> \Oc_\G(G_k)@>\Delta'_2(G_k)>> \ldots @>\Delta'_r(G_k)>>\Oc_\G(G_k)\\
@VV\xi_g^*V\ @VV\xi_g^*V @.@VV\xi_g^*V\\
\Oc_\G(U) @>\Delta'_1(U)>> \Oc_\G(U)@>\Delta'_2(U)>>\ldots @>\Delta'_r(U)>>\Oc_\G(U).
\end{CD}\end{equation}
Let   $\Delta\in\U_T$ and $\alpha\in\Ac_R$  be such that $(g,1)\in U^{(R)}_\alpha$. Take $P\in\field[x_0,\ldots,x_N]$  homogeneous.  Then (\ref{E17}) leads to 
\begin{equation}
\label{E18}
\Dc(\Delta)\Bigg(\frac{P\big(T^{(R)}_\alpha(\phi_{k}(g),\y)\big)}{T^{(R)}_{k,\alpha}(\phi_{k}(g),\y)^{\deg(P)}}\Bigg)\Bigg|_{\y=\phi_0(1)}=\Bc_k(\Delta)\bigg(\frac{P(\x)}{x_k^{\deg(P)}}\bigg)\Bigg|_{\x=\phi_k(g)}.
\end{equation}
Set $w^{(T)}_{k,\alpha}(P):=T^{(R)}_{k,\alpha}(\x,\phi_0(1))^{\deg(P)+T}$ and $v^{(T)}_{k}(P):=x_k^{\deg(P)+T(c_6-1)}$.  We deduce from (\ref{E18}) that for all $\Delta\in\U_T$, $\alpha\in\Ac_R$ and $g\in G$ such that $(g,1)\in U_\alpha^{(R)}$
\begin{align*}
w^{(T)}_{k,\alpha}(P)\cdot v^{(T)}_{k}(P)\cdot \Dc(\Delta)\Bigg(\frac{P\big(T^{(R)}_\alpha(\phi_{k}(g),\y)\big)}{T^{(R)}_{k,\alpha}(\phi_{k}(g),\y)^{\deg(P)}}\Bigg)&\\
-w^{(T)}_{k,\alpha}(P)\cdot v^{(T)}_{k}(P)\cdot \Bc_k(\Delta)\bigg(\frac{P(\x)}{x_k^{\deg(P)}}\bigg)\Bigg|_{(\x,\y)=(\phi_k(g),\phi_0(1))}&=0;
\end{align*}
since $U_\alpha^{(R)}$ is dense in $\G\times G$, we conclude that
\begin{equation}
\label{E19}
 v^{(T)}_k(P)\cdot P^{(k)}_{\Delta,\alpha}-w^{(T)}_{k,\alpha}(P)\cdot P^{\Delta}_k\in I(\G).
\end{equation}
Now we start the induction.  Lemma \ref{R14} iv) yields the equality if $T=0$. Thus from know on we assume that $T\in\Nat$ and $\Cb^{T-1}(I)=\Db^{T-1}(I)$. With the notation as above, define
\begin{align*}
K_T&:=\Ib\Big(w^{(T)}_{k,\alpha}(P)\cdot P^{\Delta'}_{k},I(\G):\; P\in \Ih,\;k\in\{0,\ldots, N\},\,\alpha\in\Ac_R,\;\Delta'\in \U_T\Big)\\
 I_T&:=\Ib\Big(w^{(T)}_{k,\alpha}(P)\cdot P^{\Delta'}_{k}, P^{\Delta}_{k},I(\G):\; P\in \Ih,\;k\in\{0,\ldots, N\},\\
 &\qquad\qquad \alpha\in\Ac_R,\; \Delta\in\U(\bfr,T-1),\;\Delta'\in \U_T\Big);
\end{align*}
hence it is clear that $I_T\subseteq \Cb^T(I)$. See that Lemma \ref{R16} i) gives
 \begin{align}
 \label{E20}
 I_T&\supseteq\Ib\Big(w^{(T)}_{k,\alpha}(P)\cdot x_l^{T(c_6-1)}\cdot P^{\Delta'}_{k}, P^{\Delta}_{k},I(\G):\; P\in \Ih,\;l,k\in\{0,\ldots, N\},\nonumber\\
 &\qquad\qquad \alpha\in\Ac_R,\; \Delta\in\U(\bfr,T-1),\;\Delta'\in \U_T\Big)\nonumber\\
 &=\Ib\Big(w^{(T)}_{k,\alpha}(P)\cdot x_k^{T(c_6-1)}\cdot P^{\Delta'}_{l}, P^{\Delta}_{k},I(\G):\; P\in \Ih,\;l,k\in\{0,\ldots, N\},\nonumber\\
 &\qquad\qquad \alpha\in\Ac_R,\;\Delta\in\U(\bfr,T-1),\;\Delta'\in \U_T\Big).
\end{align}
Now we have that 
\begin{equation*}
\Z\Big(w^{(T)}_{k,\alpha}(P)\cdot x_k^{T(c_6-1)},\;I(\G):\;k\in\{0,\ldots, N\},\;\alpha\in\Ac_R\Big)=\emptyset.
\end{equation*}
Hence Proposition \ref{R11} and (\ref{E20}) yield 
\begin{align*}
 I_T&\supseteq\Ib\Big(P^{\Delta'}_{k}, P^{\Delta}_{k},I(\G):\; P\in \Ih,\;k\in\{0,\ldots, N\},\\
 &\qquad\qquad \alpha\in\Ac_R,\; \Delta\in\U(\bfr,T-1),\;\Delta'\in \U_T\Big)
\end{align*}
and consequently we have the equality. Then, inasmuch as   $\U(\bfr,T-1)\cup \U_T$ generates $\U(\bfr,T)$, we get that $I_T=\Cb^T(I)$.  Name
\begin{align*}
M_T&:=\Ib\Big(v^{(T)}_{k}(P)\cdot P^{(k)}_{\Delta',\alpha},I(\G):\; P\in \Ih,\;k\in\{0,\ldots,N\},\,\alpha\in\Ac_R,\;\Delta'\in \U_T\Big)\\
 J_T&:=\Ib\Big(v^{(T)}_{k}(P)\cdot P^{(k)}_{\Delta',\alpha}, P^{(k)}_{\Delta,\alpha},I(\G):\; P\in \Ih,\;k\in\{0,\ldots,N\},\\
 &\qquad\qquad\alpha\in\Ac_R,\; \Delta\in\U(\bfr,T-1),\;\Delta'\in \U_T\Big)
\end{align*}
 and thus $J_T\subseteq \Db^T(I)$. Lemma \ref{R16} ii)  leads to 
 \begin{align}
 \label{E21}
 J_T&\supseteq \Ib\Big(v^{(T)}_{k}(P)\cdot T^{(R)}_{l,\alpha}(\x,\phi_0(1))^{T}\cdot P^{(k)}_{\Delta',\alpha}, P^{(k)}_{\Delta,\alpha},I(\G):\; P\in \Ih,\nonumber\\
 &\qquad\qquad l, k\in\{0,\ldots,N\},\; \alpha\in\Ac_R,\; \Delta\in\U(\bfr,T-1),\;\Delta'\in \U_T\Big)\nonumber\\
 &=\Ib\Big(v^{(T)}_{k}(P)\cdot T^{(R)}_{k,\alpha}(\x,\phi_0(1))^{T}\cdot P^{(l)}_{\Delta',\alpha}, P^{(k)}_{\Delta,\alpha},I(\G):\; P\in \Ih,\nonumber\\
 &\qquad\qquad l, k\in\{0,\ldots,N\},\; \alpha\in\Ac_R,\; \Delta\in\U(\bfr,T-1),\;\Delta'\in \U_T\Big).
\end{align}
Since \begin{equation*}
\Z\Big(v^{(T)}_{k}(P)\cdot T^{(R)}_{k,\alpha}(\x,\phi_0(1))^{T},\;I(\G):\;k\in\{0,\ldots,N\},\;\alpha\in\Ac_R\Big)=\emptyset,
\end{equation*}
 Proposition \ref{R11} and (\ref{E21}) yield
 \begin{align*} 
  J_T&\supseteq\Ib\Big(P^{(k)}_{\Delta',\alpha}, P^{(k)}_{\Delta,\alpha},I(\G):\; P\in \Ih,\;k\in\{0,\ldots,N\},\\
 &\qquad\qquad\alpha\in\Ac_R,\; \Delta\in\U(\bfr,T-1),\;\Delta'\in \U_T\Big)
 \end{align*}and trivially we have the equality. Moreover, since $\U(\bfr,T-1)\cup \U_T$ generates $\U(\bfr,T)$, we conclude that  that $J_T=\Db^T(I)$. From Corollary \ref{R10} ii)
 \begin{equation}
 \label{E57}
 I_T=(K_T,\Cb^{T-1}(I))\qquad\text{and}\qquad J_T=(M_T,\Db^{T-1}(I)).
 \end{equation}
 See that (\ref{E19}) yields
 \begin{equation}
 \label{E58}
 K_T=M_T
 \end{equation}
  From the equality $\Cb^{T-1}(I)=\Db^{T-1}(I)$, (\ref{E57}) and (\ref{E58}), we conclude that  $I_T=J_T$  and this completes the induction.
\end{proof}
\begin{lem}
\label{R19}
Let $I$ be a homogeneous ideal of $\field[x_0,\ldots,x_N]$ and  $T\in \Nat\cup\{0\}$. Then
\begin{equation*}
\Db^{T}(I)=\Eb^{T}(I).
\end{equation*}
\end{lem}
\begin{proof}
The proof shall be done by induction on $T$. Trivially $\Db^{0}(I)=\Eb^{0}(I)$; hence we assume form now on that $T\in\Nat$ and $\Db^{T-1}(I)=\Eb^{T-1}(I)$. Let $P\in\field[x_0,\ldots,x_N]$ be homogeneous, $k\in\{0,\ldots,N\}$, $\alpha\in \Ac_R$ and $\Delta\in \U_T$. A straight consequence of applying the Leibniz's rule to the  the product $T^{(R)}_{k,\alpha}(\x,\y)^{\deg(P)+T}\cdot\frac{P(T^{(R)}_\alpha(\x,\y))}{T_{k,\alpha}(\x,\y)^{\deg(P)}}$
 is  that
\begin{equation}
\label{E22}
T^{(R)}_{k,\alpha}(\x,\phi_0(1))^{T}\cdot P_{\Delta,\alpha}-P^{(k)}_{\Delta,\alpha}\in \Db^{T-1}(P).
\end{equation}
Set 
\begin{align*}
K_T&:=\Ib\Big(T^{(R)}_{k,\alpha}(\x,\phi_0(1))^{T}\cdot P_{\Delta',\alpha},I(\G):\; P\in \Ih,\;k\in\{0,\ldots,N\},\\
&\qquad\qquad \alpha\in\Ac_R,\; \Delta'\in \U_T\Big)\\
 I_T&:=\Ib\Big(T^{(R)}_{k,\alpha}(\x,\phi_0(1))^{T}\cdot P_{\Delta',\alpha}, P_{\Delta,\alpha},I(\G):\; P\in \Ih,\;k\in\{0,\ldots,N\},\\
 &\qquad\qquad \alpha\in\Ac_R,\; \Delta\in\U(\bfr,T-1),\;\Delta'\in \U_T\Big).
\end{align*}
By Corollary \ref{R10} ii)
\begin{equation}
\label{E59}
I_T=(K_T,\Eb^{T-1}(I))
\end{equation}
Since  $\Db^{T-1}(I)=\Eb^{T-1}(I)$ and $\U(\bfr,T-1)\cup \U_T$ generates $\U(\bfr,T)$, we get that  (\ref{E22}) and (\ref{E59})  imply $\Db^T(I)=I_T$. Define
\begin{align*}
 J_T&:=\Ib\Big(T^{(R)}_{k,\beta}(\x,\phi_0(1))^{\deg(P)}\cdot T^{(R)}_{k,\alpha}(\x,\phi_0(1))^{T}\cdot P_{\Delta',\alpha}, P_{\Delta,\alpha},I(\G):\; P\in \Ih,\\
 &\qquad\qquad k\in\{0,\ldots,N\},\;\alpha,\beta\in\Ac_R,\; \Delta\in\U(\bfr,T-1),\;\Delta'\in \U_T\Big).
\end{align*}
Note that $J_T\subseteq I_T\subseteq \Eb^T(I)$ and therefore  it is enough to show that $J_T=\Eb^T(I)$ to conclude the proof.  A consequence of Lemma \ref{R16} iii) is that 
\begin{align}
\label{E23}
 J_T&=\Ib\Big( T^{(R)}_{k,\alpha}(\x,\phi_0(1))^{T+\deg(P)}\cdot P_{\Delta',\beta}, P_{\Delta,\alpha},I(\G):\; P\in \Ih,\;k\in\{0,\ldots,N\},\nonumber\\
 &\qquad\qquad \alpha,\beta\in\Ac_R,\; \Delta\in\U(\bfr,T-1),\;\Delta'\in \U_T\Big).
\end{align}
Since
\begin{equation*}
\Z\Big(T^{(R)}_{k,\alpha}(\x,\phi_0(1))^{T+\deg(P)},\;I(\G):\;k\in\{0,\ldots, N\},\;\alpha\in\Ac_R\Big)=\emptyset,
\end{equation*}
Proposition \ref{R11} and (\ref{E23}) lead to the equality
\begin{align*}
 J_T&=\Ib\Big( P_{\Delta',\beta}, P_{\Delta,\alpha},I(\G):\; P\in \Ih,\;k\in\{0,\ldots,N\},\nonumber\\
 &\qquad\qquad \alpha,\beta\in\Ac_R,\; \Delta\in\U(\bfr,T-1),\;\Delta'\in \U_T\Big).
\end{align*}
  and this implies $J_T=\Eb^T(I)$ since $\U(\bfr,T-1)\cup \U_T$ generates $\U(\bfr,T)$.
\end{proof}

\begin{lem}
 \label{R20}
 Let $I$ be a homogeneous ideal of $\field[x_0,\ldots,x_N]$ and $T\in\Nat\cup\{0\}$. Then
  \begin{equation*}
 \Eb^{T}(I)=\Eb^{T}(\Ib(I)).
 \end{equation*}
 \end{lem}
 \begin{proof}
 The inclusion $\Eb^T(I)\subseteq \Eb^{T}(\Ib(I))$ is trivial since $I\subseteq \Ib(I)$. For a  homogeneous ideal $K$ of $\field[x_0,\ldots,x_N]$ and $t\in\Nat\cup\{0\}$, write 
 \begin{equation*}
 K_t:=\big(P_{\Delta,\alpha}, I(G):\;P\in\,^hK,\,\alpha\in \Ac_R,\,\Delta\in \U(\bfr,t)\big).
 \end{equation*}
 By Corollary \ref{R10} and the Primary Decomposition Theorem, there is  a homogeneous ideal $J$ of $\field[x_0,\ldots,x_N]$ such that $\Ib(I)\cap J=I$ and $\Z(J)=\emptyset$. We need to show that for all $r,s\in\Nat\cup\{0\}$ with $s\leq r$ we have
 \begin{equation}
 \label{E24}
 \Ib\big(\Ib(I)_{r-s}\cdot J_s\big)\subseteq \Ib\big((\Ib(I)\cdot J)_r\big).
 \end{equation}
 We prove this equation by induction on $r$. Trivially (\ref{E24}) holds true for $r=0$. Now  assume that $r\in\Nat$ and that the result holds true for $r-1$. If $s>0$, then 
 \begin{equation*}
 \Z(J_s)=\Z(J_{s-1})=\emptyset
 \end{equation*}
 inasmuch as 
  \begin{equation*}
\Ib(J)\subseteq \Ib(J_{s-1})\subseteq \Ib(J_s).
 \end{equation*}
 Then Proposition \ref{R11} yields
 \begin{equation*}
\Ib\big(\Ib(I)_{r-s}\cdot J_s\big)=\Ib\big(\Ib(I)_{r-s}\big)=\Ib\big(\Ib(I)_{r-s}\cdot J_{s-1}\big).
 \end{equation*}
 and the hypothesis of induction implies
 \begin{align}
  \label{E25}
 \Ib\big(\Ib(I)_{r-s}\cdot J_s\big)&=\Ib\big(\Ib(I)_{r-s}\cdot J_{s-1}\big)\nonumber\\
 &\subseteq\Ib\big((\Ib(I)\cdot J)_{r-1}\big)\nonumber\\
 &\subseteq\Ib\big((\Ib(I)\cdot J)_{r}\big).
 \end{align}
 It remains to show (\ref{E24}) when $s=0$. For all $t\in\{0,\ldots,r\}$, if $P\in \Ib(I)$, $Q\in J$, $\alpha\in\Ac_R$ and $\Delta\in\U_t$, the Leibniz rule and (\ref{E25}) lead to
 \begin{equation*}
  P_{\Delta,\alpha}\cdot Q^{(R)}_{1,\alpha}-(P\cdot Q)_{\Delta,\alpha}\in \Ib\big((\Ib(I)\cdot J)_r\big).
 \end{equation*}
 Thus, since $(P\cdot Q)_{\Delta,\alpha}\in \Ib\big((\Ib(I)\cdot J)_r\big)$, we deduce that
   \begin{equation}
 \label{E60}
 P_{\Delta,\alpha}\cdot Q^{(R)}_{1,\alpha}\in \Ib\big((\Ib(I)\cdot J)_r\big).
 \end{equation}
Inasmuch as $\bigcup_{t=0}^r\U_t$ generates $\U(\bfr,r)$, we conclude from (\ref{E60}) that 
 \begin{equation*}
  \Ib\big(\Ib(I)_{r}\cdot J_0\big)\subseteq \Ib\big((\Ib(I)\cdot J)_r\big)
 \end{equation*}
 which completes the proof by induction of (\ref{E24}). We come back to the proof and we finish it as follows
 \begin{align*}
 \Eb^T(\Ib(I))&=\Ib\big(\Ib(I)_T\big)\\
 &=\Ib\big(\Ib(I)_T\cdot J\big)&\text{by Proposition \ref{R11}}\\
 &\subseteq\Ib\big((\Ib(I)\cdot J)_T\big)&\text{by (\ref{E24})}\\
 &\subseteq\Ib\big(I_T\big)\\
 &=\Eb^T(I).
 \end{align*}
  \end{proof}
All the effort we spent showing that $\Bb^T(I)=\Cb^T(I)=\Db^T(I)=\Eb^T(I)$ is rewarded with the following claim.
\begin{cor}
\label{R21}
Let $I$ be a homogeneous ideal and $T,T'\in\Nat\cup\{0\}$. Then 
\begin{equation*}
\Eb^{T+T'}(I)=\Eb^T\big(\Eb^{T'}(I)\big).
\end{equation*}
\end{cor}
\begin{proof}
For any   $t\in\Nat\cup\{0\}$, set
\begin{align*}
I_t:=&\Big(P^{(\Delta'_1,\ldots,\Delta'_r)}_{(k_1,\ldots,k_r)},I,I(\G):\; P\in \Ih,\;r\in\{1,\ldots, t\},\\
&\qquad  k_1\ldots,k_r\in \{0,\ldots,N\},\; \Delta'_1,\ldots,\Delta'_r\in\bfr\Big)
\end{align*}
By Lemma \ref{R17}, Lemma \ref{R18} and Lemma \ref{R19}, we get that for all $t\in\Nat\cup\{0\}$
\begin{equation}
\label{E26}
\Bb^t(I)=\Eb^t(I).
\end{equation}
Then
\begin{align*}
\Eb^{T}(\Eb^{T'}(I))&=\Eb^T(\Bb^{T'}(I))&\text{by (\ref{E26})}\\
&=\Eb^T(I_{T'})&\text{by Lemma \ref{R20}}\\
&=\Bb^T(I_{T'})&\text{by (\ref{E26})}\\
&=\Bb^{T+T'}(I)\\
&=\Eb^{T+T'}(I)&\text{by (\ref{E26})}.
\end{align*}
\end{proof}
\begin{prop}
\label{R22}
Let $I$ be a homogeneous ideal of $\field[x_0,\ldots,x_N]$,  $T\in\Nat\cup\{0\}$ and $g\in G$. Then 
\begin{equation*}
\Tb_{L_g}\big(\Eb^T(I)\big)=\Ib\big(\partial_{L_g}^T(I)\big)=\Eb^T\big(\Tb_{L_g}(I)\big).
\end{equation*}
\end{prop}
\begin{proof}
Recall that if $z\in \G$, then $k_z:=\min\{k\in\{0,\ldots,N\}:\;z\in G_k\}$. Take $P\in \field[x_0,\ldots,x_N]$ homogeneous, $\alpha,\gamma\in \Ac_R$, $\beta, \delta\in \Ac_L$ and $k\in\{0,\ldots,N\}$, and write
\begin{align*} 
  v^{(\delta,\gamma)}_{k}(P)&:=T^{(L)}_{k,\delta}\Big(\phi_{k_g}(g),T^{(R)}_\gamma(\x,\y)\Big)^{\deg(P)}\\
  w^{(\alpha,\beta)}_{k}(P)&:=T^{(R)}_{k,\alpha}\Big(T^{(L)}_\beta(\phi_{k_g}(g),\x),\y\Big)^{\deg(P)}\\
  Q^{(\delta,\gamma)}(P)&:=P\Big(T^{(L)}_\delta\Big(\phi_{k_g}(g),T^{(R)}_\gamma(\x,\y)\Big)\Big)\\
  R^{(\alpha,\beta)}(P)&:=P\Big(T^{(R)}_\alpha\Big(T^{(L)}_\beta(\phi_{k_g}(g),\x),\y\Big)\Big).
   \end{align*}
Note that for all $\Delta\in \U(\bfr)$
\begin{align*}
 \Dc(\Delta)\big(Q^{(\delta,\gamma)}(P)\big)\Big|_{\y=\phi_0(1)}=\big(P^{(L)}_{g,\delta}\big)_{\Delta,\gamma}&\qquad  v^{(\delta,\gamma)}_{k}(P)\Big|_{\y=\phi_0(1)}=\Big(\big(x_k^{\deg(P)}\big)^{(L)}_{g,\delta}\Big)^{(R)}_{1,\gamma}\\
 \Dc(\Delta)\big(R^{(\alpha,\beta)}(P)\big)\Big|_{\y=\phi_0(1)}=\big(P_{\Delta,\alpha}\big)^{(L)}_{g,\beta}&\qquad  w^{(\alpha,\beta)}_{k}(P)\Big|_{\y=\phi_0(1)}=\Big(\big(x_k^{\deg(P)}\big)^{(R)}_{1,\alpha}\Big)^{(L)}_{g,\beta}.
\end{align*}
For $t\in \Nat\cup\{0\}$ define the ideal of $\field[x_0,\ldots,x_N]$
\begin{align*}
I_{t}:=&\bigg(\big(P^{(L)}_{g,\delta}\big)_{\Delta,\gamma},I(\G):\;P\in \Ih,\;\gamma\in\Ac_R,\;\delta\in \Ac_L,\;\Delta\in \U(\bfr,t)\bigg).
\end{align*}
We shall prove the following equality by induction on $T$
\begin{equation}
\label{E27}
\Ib\big(\partial^T_{L_g}(I)\big)=\Ib(I_{T}).
\end{equation}
 From Lemma \ref{R14} iii) we get that (\ref{E27}) holds true for $T=0$; thus we may assume that $T\in \Nat$ and  $\Ib\big(\partial^{T-1}_{L_g}(I)\big)=\Ib(I_{T-1})$. Take $(z,w)\in \G\times G$ and call $ u:=(\phi_{k_z}(z),\phi_{k_w}(w))$.  First assume that 
 \begin{equation}
 \label{E28}
 (g,z)\in U^{(L)}_\beta,\quad(gz,w)\in U^{(R)}_\alpha,\quad(z,w)\in U^{(R)}_\gamma,\quad (g,zw)\in U^{(L)}_\delta;
 \end{equation} 
  then
 \begin{align}
 \label{E29}
 v^{(\delta,\gamma)}_{k}(P)\cdot R^{(\alpha,\beta)}(P)-
w^{(\alpha,\beta)}_{k}(P)\cdot   Q^{(\delta,\gamma)}(P)\bigg|_{(\x,\y)=u}=0
 \end{align}
 since 
 \begin{align*}
 \phi\big(g(zw)\big)&=\bigg[ T^{(L)}_{0,\delta}\Big(\phi_{k_g}(g),T^{(R)}_\gamma(\x,\y)\Big):\ldots: T^{(L)}_{N,\delta}\Big(\phi_{k_g}(g),T^{(R)}_\gamma(\x,\y)\Big)\bigg]\bigg|_{(\x,\y)=u}\\
 &=\bigg[T^{(R)}_{0,\alpha}\Big(T^{(L)}_\beta(\phi_{k_g}(g),\x),\y\Big):\ldots:T^{(R)}_{N,\alpha}\Big(T^{(L)}_\beta(\phi_{k_g}(g),\x),\y\Big)\bigg]\bigg|_{(\x,\y)=u}\\
 &=\phi\big((gz)w\big).
 \end{align*}
 Now see that if  (\ref{E28}) is not satisfied, then  (\ref{E29}) still holds true inasmuch as both addends of (\ref{E29}) would be zero. Thus, since (\ref{E29}) holds true for all $(z,w)\in \G\times G$, we obtain that
 \begin{equation*}
  v^{(\delta,\gamma)}_{k}(P)\cdot R^{(\alpha,\beta)}(P)-
w^{(\alpha,\beta)}_{k}(P)\cdot   Q^{(\delta,\gamma)}(P)\bigg|_{\y=\phi_0(1)}\in I(\G)
 \end{equation*}
 and therefore we get that for all $\Delta\in \U(\bfr)$
 \begin{equation}
 \label{E30}
  \Dc(\Delta)\Big(v^{(\delta,\gamma)}_{k}(P)\cdot R^{(\alpha,\beta)}(P)-
w^{(\alpha,\beta)}_{k}(P)\cdot   Q^{(\delta,\gamma)}(P)\Big)\bigg|_{\y=\phi_0(1)}\in I(\G)
 \end{equation}
  Let $\Delta$ be in $\U(\bfr,T)$ and note that, by the Leibniz's rule, there exist $R\in \partial^{T-1}_{L_g}(I)$  and $Q\in I_{T-1}$ such that  
 \begin{align}
 \label{E31}
 \Dc(\Delta)\Big(v^{(\delta,\gamma)}_{k}(P)\cdot R^{(\alpha,\beta)}(P)\Big)\bigg|_{\y=\phi_{0}(1)}- \Big(\big(x_k^{\deg(P)}\big)^{(L)}_{g,\delta}\Big)^{(R)}_{1,\gamma}\cdot \big(P_{\Delta,\alpha}\big)^{(L)}_{g,\beta}-R\in I(\G)\nonumber\\
  \Dc(\Delta)\Big(w^{(\alpha,\beta)}_{k}(P)\cdot   Q^{(\delta,\gamma)}(P)\Big)\bigg|_{\y=\phi_{0}(1)}-  \Big(\big(x_k^{\deg(P)}\big)^{(R)}_{1,\alpha}\Big)^{(L)}_{g,\beta}\cdot\big(P^{(L)}_{g,\delta}\big)_{\Delta,\gamma} -Q\in I(\G).
 \end{align} 
Define
 \begin{align}
 J_{T}:=& \bigg(\Dc(\Delta)\Big(v^{(\delta,\gamma)}_{k}(P)\cdot R^{(\alpha,\beta)}(P)\Big)\Big|_{\y=\phi_{0}(1)},\partial^{T-1}_{L_g}(I):\; P\in \Ih,\nonumber\\
&\quad  \alpha,\gamma\in\Ac_R,\;\beta,\delta\in \Ac_L,\;k\in\{0,\ldots,N\},\;\Delta\in \U(\bfr,T)\bigg)\nonumber\\
K_{T}:=&\bigg( \Dc(\Delta)\Big(w^{(\alpha,\beta)}_{k}(P)\cdot   Q^{(\delta,\gamma)}(P)\Big)\Big|_{\y=\phi_{0}(1)},I_{T-1}:\; P\in \Ih,\nonumber\\
&\quad \alpha,\gamma\in\Ac_R,\;\beta,\delta\in \Ac_L,\;k\in\{0,\ldots,N\},\;\Delta\in \U(\bfr,T)\bigg)\nonumber.
 \end{align}
 Thus (\ref{E31}) leads to 
  \begin{align}
  \label{E32}
 J_{T}=& \bigg(\Big(\big(x_k^{\deg(P)}\big)^{(L)}_{g,\delta}\Big)^{(R)}_{1,\gamma}\cdot \big(P_{\Delta,\alpha}\big)^{(L)}_{g,\beta},\;\partial^{T-1}_{L_g}(I):\; P\in \Ih,\;\alpha,\gamma\in\Ac_R,\nonumber\\
&\qquad  \beta,\delta\in \Ac_L,\;k\in\{0,\ldots,N\},\;\Delta\in \U(\bfr,T)\bigg)\nonumber\\
K_{T}=&\bigg( \Big(\big(x_k^{\deg(P)}\big)^{(R)}_{1,\alpha}\Big)^{(L)}_{g,\beta}\cdot \big(P^{(L)}_{g,\delta}\big)_{\Delta,\gamma},I_{T-1}:\; P\in \Ih,\;\alpha,\gamma\in\Ac_R,\nonumber\\
&\qquad \beta,\delta\in \Ac_L,\;k\in\{0,\ldots,N\},\;\Delta\in \U(\bfr,T)\bigg).
 \end{align}
 On one hand,  it is easily seen that
 \begin{align*}
 \Z\bigg(\Big(\big(x_k^{\deg(P)}\big)^{(L)}_{g,\delta}\Big)^{(R)}_{1,\gamma},\;I(\G):\;\gamma\in\Ac_R,\;\delta\in\Ac_L,\;k\in\{0,\ldots,N\}\bigg)&=\\
 \Z\bigg(\Big(\big(x_k^{\deg(P)}\big)^{(R)}_{1,\alpha}\Big)^{(L)}_{g,\beta},\;I(\G):\;\alpha\in\Ac_R,\;\beta\in\Ac_L,\;k\in\{0,\ldots,N\}\bigg)&=\emptyset,
 \end{align*}
 and hence  Proposition \ref{R11} and (\ref{E32}) yield
 \begin{equation}
 \label{E33}
  \Ib(\partial^{T}_{L_g}(I))=\Ib( J_{T})\qquad\text{and}\qquad\Ib(I_{T})=\Ib(K_{T}).
 \end{equation}
  On the other hand, (\ref{E30}) implies that for all $P\in \Ih$, 
$\alpha,\gamma\in\Ac_R$, $\beta,\delta\in \Ac_L$, $k\in\{0,\ldots,N\}$ and $\Delta\in \U(\bfr,T)$, the two polynomials $\Dc(\Delta)\big(v^{(\delta,\gamma)}_{k}(P)\cdot R^{(\alpha,\beta)}(P)\big)\big|_{\y=\phi_{0}(1)}$ and $\Dc(\Delta)\big(w^{(\alpha,\beta)}_{k}(P)\cdot   Q^{(\delta,\gamma)}(P)\big)\big|_{\y=\phi_{0}(1)}$ are equal up to an element of $I(\G)$; then (\ref{E31}) and the hypothesis of induction give
 \begin{align}
 \label{E34}
 \Big(\big(x_k^{\deg(P)}\big)^{(L)}_{g,\delta}\Big)^{(R)}_{1,\gamma}\cdot \big(P_{\Delta,\alpha}\big)^{(L)}_{g,\beta}\in K_{T}\qquad\text{and}\qquad \Big(\big(x_k^{\deg(P)}\big)^{(R)}_{1,\alpha}\Big)^{(L)}_{g,\beta}\cdot\big(P^{(L)}_{g,\delta}\big)_{\Delta,\gamma}\in J_{T}.
 \end{align}
 See that (\ref{E34}) and the hypothesis of induction let us assert, using Corollary \ref{R10} ii), that
 $\Ib( J_{T})=\Ib( K_{T})$, and this equality jointly with (\ref{E33}) concludes the proof of (\ref{E27}). Now we complete the proof the statement
 \begin{align*}
 \Tb_{L_g}\big(\Eb^T(I)\big)&= \Ib(\partial^{T}_{L_g}(I))&\text{by Lemma \ref{R13}}\\
 &=\Ib(I_{T})&\text{by (\ref{E27})}\\
 &=\Eb^T\big(\Tb_{L_g}(I)\big)&\text{by Lemma \ref{R20}}.
 \end{align*} 
\end{proof}
Before we proceed with the analogous result  to Proposition \ref{R22} for $\partial_{R_g}^T(I)$, we need a  technical result.
\begin{lem}
\label{R23}
Let $P\in \field[x_0,\ldots,x_N]$ be homogeneous,  $\alpha,\beta,\gamma,\delta\in \Ac_R$, $k,l\in\{0,\ldots,N\}$ and $\Delta\in\U(\bfr,T)\setminus\U(\bfr,T-1)$. For all  $g\in G$  such that $\Ad(g)(\bfr)\subseteq \bfr$
\begin{equation*}
\Big(\big(x_k^{\deg(P)+T}\big)^{(R)}_{1,\alpha}\Big)^{(R)}_{g,\beta}\cdot \big(P_{g,\delta}^{(R)}\big)^{(l)}_{\Ad(g)(\Delta),\gamma}-\Big(\big(x_l^{\deg(P)+T}\big)^{(R)}_{g,\delta}\Big)^{(R)}_{1,\gamma}\cdot\big(P_{\Delta,\alpha}^{(k)}\big)_{g,\beta}^{(R)}\in I(\G).
\end{equation*}
\end{lem}
\begin{proof}
Write 
\begin{equation*}
  W_1:=\Big\{z\in\G:\; zg\in G_l,\;(z,1)\in U^{(R)}_{\gamma},\; (z,g)\in U^{(R)}_{\delta}\Big\}.
   \end{equation*} 
 First assume that $W_1$ is empty. In this case, for any homogeneous polynomial $Q\in\field[x_0,\ldots,x_N]$, $t\in\Nat\cup\{0\}$ and   $\Delta'\in \U(\bfr,t)\setminus\U(\bfr,t-1)$, the polynomial
\begin{align*}
 \big(Q(\x)_{g,\delta}^{(R)}\big)^{(l)}_{\Delta',\gamma}&=\\
  T^{(R)}_{l,\delta}\Big(T^{(R)}_\gamma(\x,\y),\phi_{k_g}(g)\Big)^{\deg(Q)+t}\cdot\Dc(\Delta')\Bigg(\frac{Q\Big(T^{(R)}_{\delta}\Big(T^{(R)}_\gamma(\x,\y),\phi_{k_g}(g)\Big)\Big)}{T^{(R)}_{l,\delta}\Big(T^{(R)}_\gamma(\x,\y),\phi_{k_g}(g)\Big)^{\deg(Q)}}\Bigg)\Bigg|_{\y=\phi_{0}(1)}&
\end{align*}evaluated in any $z\in\G$ is zero and therefore it is in $I(\G)$. Hence the statement of the lemma  is true in this case since $ \Big(\big(x_l^{\deg(P)+T}\big)^{(R)}_{g,\delta}\Big)^{(R)}_{1,\gamma}, \big(P_{g,\delta}^{(R)}\big)^{(l)}_{\Ad(g)(\Delta)}\in I(G)$. Now write
\begin{equation*}
  W_2:=\Big\{z\in\G:\; zg\in G_k,\;(zg,1)\in U^{(R)}_{\gamma},\; (z,g)\in U^{(R)}_{\delta}\Big\}.
   \end{equation*}
  If $W_2$ is empty,  we deduce proceeding as above that $\Big(\big(x_k^{\deg(P)+T}\big)^{(R)}_{1,\alpha}\Big)^{(R)}_{g,\beta}, \big(P_{\Delta,\alpha}^{(k)}\big)_{g,\beta}^{(R)}\in I(\G)$ and the statement is true also in this case. From now on we assume that  $W_1$ and $W_2$ are not empty. Take $\alpha',\gamma',\delta'\in \Ac_R$ and $\beta'\in\Ac_L$ such that 
   \begin{equation*}
  W_3:=\Big\{z\in\G:\; zg\in G_l,\;(1,g)\in U^{(R)}_{\gamma'},\; (z,g)\in U^{(R)}_{\delta'}\Big\}
   \end{equation*} 
and 
 \begin{equation*}
  W_4:=\Big\{z\in\G:\; zg\in G_k,\;(g,1)\in U^{(L)}_{\beta'},\quad (z,g)\in U^{(R)}_{\alpha'}\Big\}
     \end{equation*} 
     are not empty. Thus  $U:=\bigcap_{i=1}^4W_i$ is an open dense subset of $\G$. For all  $z\in U$, we get that 
    \begin{align*}
       \frac{P\Big(T^{(R)}_{\delta}\Big(T^{(R)}_\gamma(\phi_{k_z}(z),\x),\phi_{k_g}(g)\Big)\Big)}{T^{(R)}_{l,\delta}\Big(T^{(R)}_\gamma(\phi_{k_z}(z),\x),\phi_{k_g}(g)\Big)^{\deg(P)}}= \frac{P\Big(T^{(R)}_{\delta'}\Big(\phi_{k_z}(z),T^{(R)}_{\gamma'}(\x,\phi_{k_g}(g))\Big)}{T^{(R)}_{l,\delta'}\Big(\phi_{k_z}(z),T^{(R)}_{\gamma'}(\x,\phi_{k_g}(g))\Big)^{\deg(P)}}
   \end{align*}
    as a regular function in a neighbourhood of $1$; hence
     \begin{align}
    \label{E35}
    \Dc(\Ad(g)(\Delta))\Bigg(\frac{P\Big(T^{(R)}_{\delta}\Big(T^{(R)}_\gamma(\x,\y),\phi_{k_g}(g)\Big)\Big)}{T^{(R)}_{l,\delta}\Big(T^{(R)}_\gamma(\x,\y),\phi_{k_g}(g)\Big)^{\deg(P)}}\Bigg)\Bigg|_{(\x,\y)=(\phi_{k_z}(z),\phi_{0}(1))}&=\nonumber\\
     \Dc(\Ad(g)(\Delta))\Bigg(\frac{P\Big(T^{(R)}_{\delta'}\Big(\x,T^{(R)}_{\gamma'}(\y,\phi_{k_g}(g))\Big)}{T^{(R)}_{l,\delta'}\Big(\x,T^{(R)}_{\gamma'}(\y,\phi_{k_g}(g))\Big)^{\deg(P)}}\Bigg)\Bigg|_{(\x,\y)=(\phi_{k_z}(z),\phi_{0}(1))}.&
   \end{align}
   In the same way we deduce that  for all $z\in U$
   \begin{align}
    \label{E36}
    \Dc(\Delta)\Bigg(\frac{P\Big(T^{(R)}_{\alpha}\Big(T^{(R)}_\beta(\x,\phi_{k_g}(g)),\y \Big)\Big)}{T^{(R)}_{k,\alpha}\Big(T^{(R)}_\beta(\x,\phi_{k_g}(g)),\y\Big)^{\deg(P)}}\Bigg)\Bigg|_{(\x,\y)=(\phi_{k_z}(z),\phi_{0}(1))}&=\nonumber\\
      \Dc(\Delta)\Bigg(\frac{P\Big(T^{(R)}_{\alpha'}\Big(\x,T^{(L)}_{\beta'}(\phi_{k_g}(g),\y )\Big)\Big)}{T^{(R)}_{k,\alpha'}\Big(\x,T^{(L)}_{\beta'}(\phi_{k_g}(g),\y )\Big)^{\deg(P)}}\Bigg)\Bigg|_{(\x,\y)=(\phi_{k_z}(z),\phi_{0}(1))}.&
   \end{align}
Since  $\Ad(g)(\Delta)=\eta_g^*\circ\Delta\circ\eta_{g^{-1}}^*$, we get that for all $z\in U$
 \begin{align}
    \label{E37}
    \Dc(\Ad(g)(\Delta))\Bigg(\frac{P\Big(T^{(R)}_{\delta'}\Big(\phi_{k_z}(z),T^{(R)}_{\gamma'}(\y,\x)\Big)}{T^{(R)}_{l,\delta'}\Big(\phi_{k_z}(z),T^{(R)}_{\gamma'}(\y,\x)\Big)^{\deg(P)}}\Bigg)\Bigg|_{(\x,\y)=(\phi_{k_g}(g),\phi_{0}(1))}&=\nonumber\\
    \Bc_{k_g}(\Delta)\Bigg(\frac{P\Big(T^{(R)}_{\delta'}\Big(\phi_{k_z}(z),\x\Big)}{T^{(R)}_{l,\delta'}\Big(\phi_{k_z}(z),\x\Big)^{\deg(P)}}\Bigg)\Bigg|_{\x=\phi_{k_g}(g)}.&
   \end{align}
  Since  $\Delta=\xi_g^*\circ\Delta\circ\xi_{g^{-1}}^*$, we conclude that for all  $z\in U$
   \begin{align}
    \label{E38}
           \Dc(\Delta)\Bigg(\frac{P\Big(T^{(R)}_{\alpha'}\Big(\phi_{k_z}(z),T^{(L)}_{\beta'}(\x,\y )\Big)\Big)}{T^{(R)}_{k,\alpha'}\Big(\phi_{k_z}(z),T^{(L)}_{\beta'}(\x,\y )\Big)^{\deg(P)}}\Bigg)\Bigg|_{(\x,\y)=(\phi_{k_g}(g),\phi_{0}(1))}&=\nonumber\\
          \Bc_{k_g}(\Delta)\Bigg(\frac{P\Big(T^{(R)}_{\delta'}\Big(\phi_{k_z}(z),\x\Big)}{T^{(R)}_{l,\delta'}\Big(\phi_{k_z}(z),\x\Big)^{\deg(P)}}\Bigg)\Bigg|_{\x=\phi_{k_g}(g)}.&
   \end{align}
   Set
   \begin{equation*}
v(\x,\y):= \bigg( T^{(R)}_{l,\delta}\Big(T^{(R)}_\gamma(\x,\y),\phi_{k_g}(g)\Big)\cdot T^{(R)}_{k,\alpha}\Big(T^{(R)}_\beta(\x,\phi_{k_g}(g)),\y\Big)\bigg)^{\deg(P)+T}.
\end{equation*}
By the equalities  (\ref{E35}),(\ref{E36}),(\ref{E37}) and (\ref{E38}), we conclude that for all $z\in U$
\begin{align}
 \label{E39}
v(\x,\y)\cdot \Dc(\Ad(g)(\Delta))\Bigg(\frac{P\Big(T^{(R)}_{\delta}\Big(T^{(R)}_\gamma(\x,\y),\phi_{k_g}(g)\Big)\Big)}{T^{(R)}_{l,\delta}\Big(T^{(R)}_\gamma(\x,\y),\phi_{k_g}(g)\Big)^{\deg(P)}}\Bigg)&\nonumber\\
-v(\x,\y)\cdot \Dc(\Delta)\Bigg(\frac{P\Big(T^{(R)}_{\alpha}\Big(T^{(R)}_\beta(\x,\phi_{k_g}(g)),\y \Big)\Big)}{T^{(R)}_{k,\alpha}\Big(T^{(R)}_\beta(\x,\phi_{k_g}(g)),\y\Big)^{\deg(P)}}\Bigg)&\Bigg|_{(\x,\y)=(\phi_{k_z}(z),\phi_{0}(1))}=0.
 \end{align}
 Finally $U$ is dense in $\G$ so (\ref{E39}) is true for all $z\in \G$ and the statements has been demonstrated.
\end{proof}

\begin{prop}
\label{R24}
Let $I$ be a homogeneous ideal of $\field[x_0,\ldots,x_N]$,  $T\in\Nat\cup\{0\}$ and $g\in G$. If $\Ad(g)(\bfr)\subseteq \bfr$, then 
\begin{equation*}
\Tb_{R_g}\big(\Eb^T(I)\big)=\Ib\big(\partial^T_{R_g}(I)\big)=\Eb^T\big(\Tb_{R_g}(I)\big).
\end{equation*}
\end{prop}
\begin{proof}

   For $t\in \Nat\cup\{0\}$ define the ideals 
\begin{align*}
\widehat{\partial^t_{R_g}}(I):=&\bigg(\big(P_{\Delta,\alpha}^{(k)}\big)_{g,\beta}^{(R)},I(\G):\; P\in \Ih,\;\alpha,\beta\in\Ac_R,\;k\in\{0,\ldots,N\},\;\Delta\in \U(\bfr,t)\bigg)\\
I_{t}:=&\bigg(\big(P_{g,\delta}^{(R)}\big)^{(k)}_{\Delta,\gamma},I(\G):\;P\in \Ih,\;\gamma,\delta\in\Ac_R,\;k\in\{0,\ldots,N\},\;\Delta\in \U(\bfr,t)\bigg).
\end{align*}
From Lemma \ref{R19} and Lemma \ref{R20}, it suffices to show
\begin{equation}
\label{E40}
\Tb_{R_g}\big(\Db^T(I)\big)=\Ib\Big(\widehat{\partial^T_{R_g}}(I)\Big)=\Db^T\big(\Tb_{R_g}(I)\big).
\end{equation}
The first step is to show the following equality by induction on $T$
\begin{equation}
\label{E41}
\Ib\Big(\widehat{\partial^T_{R_g}}(I)\Big)=\Ib(I_{T}).
\end{equation}
 From Lemma \ref{R13} and Lemma \ref{R14} ii),  we get that (\ref{E41}) holds  for $T=0$; thus we may assume that $T\in \Nat$ and  $\Ib\Big(\widehat{\partial^{T-1}_{R_g}}(I)\Big)=\Ib(I_{T-1})$.  Define the ideals
\begin{align*}
J_{T}:=&\bigg(\Big(\big(x_l^{\deg(P)+T}\big)^{(R)}_{g,\delta}\Big)^{(R)}_{1,\gamma}\cdot\big(P_{\Delta,\alpha}^{(k)}\big)_{g,\beta}^{(R)},\;\widehat{\partial^{T-1}_{R_g}}(I):\; P\in \Ih,\nonumber\\
&\qquad \alpha,\beta,\gamma,\delta\in\Ac_R,\;l,k\in\{0,\ldots,N\},\;\Delta\in \U(\bfr,T)\setminus\U(\bfr,T-1)\bigg)\\
K_{T}:=&\bigg(\Big(\big(x_k^{\deg(P)+T}\big)^{(R)}_{1,\alpha}\Big)^{(R)}_{g,\beta}\cdot \big(P_{g,\delta}^{(R)}\big)^{(l)}_{\Delta,\gamma},\;I_{T-1}:\; P\in \Ih,\nonumber\\
&\qquad\alpha,\beta,\gamma,\delta\in\Ac_R,\; l,k\in\{0,\ldots,N\},\;\Delta\in \U(\bfr,T)\setminus\U(\bfr,T-1)\bigg).
\end{align*}
 On one hand
 \begin{align*}
 \Z\bigg(\Big(\big(x_l^{\deg(P)+T}\big)^{(R)}_{g,\delta}\Big)^{(R)}_{1,\gamma},I(\G):\;\gamma,\delta\in\Ac_R,\;l\in\{0,\ldots,N\}\bigg)&=\\
 \Z\bigg(\Big(\big(x_k^{\deg(P)+T}\big)^{(R)}_{1,\alpha}\Big)^{(R)}_{g,\beta},I(\G):\;\alpha,\beta\in\Ac_R,\;k\in\{0,\ldots,N\}\bigg)&=\emptyset,
 \end{align*}
 and hence  Proposition \ref{R11}  yields
 \begin{equation}
 \label{E42}
  \Ib\Big(\widehat{\partial^{T}_{R_g}}(I)\Big)=\Ib( J_{T})\qquad\text{and}\qquad\Ib(I_{T})=\Ib(K_{T}).
 \end{equation}
 On the other hand, Lemma \ref{R23} asserts that  for $P\in \Ih$,\; $l,k\in\{0,\ldots,N\}$,\; $\alpha,\beta,\gamma,\delta\in\Ac_R$ and $\Delta\in \U(\bfr,T)\setminus\U(\bfr,T-1)$
 \begin{equation*}
 \Big(\big(x_l^{\deg(P)+T}\big)^{(R)}_{g,\delta}\Big)^{(R)}_{1,\gamma}\cdot\big(P_{\Delta,\alpha}^{(k)}\big)_{g,\beta}^{(R)}\in K_T\quad\text{and}\quad\Big(\big(x_k^{\deg(P)+T}\big)^{(R)}_{1,\alpha}\Big)^{(R)}_{g,\beta}\cdot \big(P_{g,\delta}^{(R)}\big)^{(l)}_{\Delta,\gamma}\in J_T.
 \end{equation*}
 Thus Corollary \ref{R10} ii) and the hypothesis of induction lead to $\Ib(J_{T})=\Ib(K_{T})$ and this equality concludes the proof of (\ref{E41}) by (\ref{E42}). We conclude the proof of (\ref{E40}) as follows 
 \begin{align*}
 \Tb_{R_g}\big(\Db^T(I)\big)&=\Ib\Big(\widehat{\partial^{T}_{R_g}}(I)\Big)&\text{by Lemma \ref{R13}}\\
 &=\Ib(I_{T})&\text{by (\ref{E41})}\\
 &=\Db^T\big(\Tb_{R_g}(I)\big)&\text{by Lemma \ref{R19} and Lemma \ref{R20}}.
 \end{align*} 
 \end{proof}
The following result is the analogous  statement to \cite[Prop. 4.3]{Philippon} that we need in the proofs of the main theorems.
\begin{cor}
\label{R25}
Let $I$ be a homogeneous ideal of $\field[x_0,\ldots,x_N]$,  $T,T'\in\Nat\cup\{0\}$ and $g,h\in G$. Then 
\begin{enumerate}
\item[i)]$\Ib\big(\partial^{T'}_{L_h}\big(\partial_{L_g}^T(I)\big)\big)=\Ib\big(\partial_{L_{gh}}^{T+T'}(I)\big)$.
\item[ii)]$\Ib\big(\partial^{T'}_{R_h}\big(\partial_{R_g}^T(I)\big)\big)=\Ib\big(\partial_{R_{hg}}^{T+T'}(I)\big)$.
\item[iii)]$\Ib\big(\partial^{T'}_{R_h}\big(\partial_{L_g}^T(I)\big)\big)=\Ib\big(\partial^{T}_{L_g}\big(\partial_{R_h}^{T'}(I)\big)\big)$.
\end{enumerate}
\end{cor}
\begin{proof}
We just show i) insomuch as the proofs of the three statements are quite similar. Then
\begin{align*}
\Ib\big(\partial^{T'}_{L_h}\big(\partial_{L_g}^T(I)\big)\big)&=\Tb_{L_h}\big(\Eb^{T'}\big(\partial_{L_g}^T(I)\big)\big)&\text{by Proposition \ref{R22}}\\
&=\Tb_{L_h}\big(\Eb^{T'}\big(\Ib\big(\partial_{L_g}^T(I)\big)\big)\big)&\text{by Lemma \ref{R20}}\\
&=\Tb_{L_h}\big(\Eb^{T'}\big(\Eb^T(\Tb_{L_g}(I))\big)\big)&\text{by Proposition \ref{R22}}\\
&=\Tb_{L_h}\big(\Eb^{T+T'}(\Tb_{L_g}(I))\big)&\text{by Corollary \ref{R21}}\\
&=\Tb_{L_h}\big(\Tb_{L_g}\big(\Eb^{T+T'}(I)\big)\big)&\text{by Proposition \ref{R22}}\\
&=\Tb_{L_{gh}}\big(\Eb^{T+T'}(I)\big)&\text{by Lemma \ref{R14} i)}\\
&=\Ib\big(\partial_{L_{gh}}^{T+T'}(I)\big)&\text{by Proposition \ref{R22}.}
\end{align*}
\end{proof}
\section{Auxiliary results}
In this section we state the main tools that will be used in the proofs of our main theorems. Denote by $B$ the connected Lie subgroup of $G(\Com)$ corresponding to the Lie subalgebra  $\bfr\otimes_\field\Com\subseteq \Lie(G(\Com))$. Remember that for an irreducible subvariety $W$ of $\G$, we define $\tau(W):=\dim(B)-\dim(W(\Com)\cap Bw)$ where  $w\in W$ is such that $W(\Com)\cap Bw$ is transverse at $w$;  if $V\subseteq \phi(\G)$ is a subvariety, $\tau(V):=\tau(\phi^{-1}(V))$. Let $Z$ be an irreducible component of a projective variety $V$ in $\Pro^N$ and   $\Oc_{Z,V}$ the local ring of $V$ along $Z$. For an homogeneous ideal $I$ of $ \field[x_0,\ldots,x_N]$ such that $Z\subseteq \Z(I)\subseteq V$, we denote by $I\Oc_{Z,V}$ its corresponding ideal in  $\Oc_{Z,V}$ and 
 \begin{equation*}
 l_{Z,V}(I):=l_{\Oc_{Z,V}}(\Oc_{Z,V}/I\Oc_{Z,V})
 \end{equation*}
 where $l_{\Oc_{Z,V}}(\Oc_{Z,V}/I\Oc_{Z,V})$ is the length of the $\Oc_{Z,V}-$module $\Oc_{Z,V}/I\Oc_{Z,V}$.  We will need the following version of Bezout's Theorem.
 \begin{thm}
 \label{R26}
 Let $F\subseteq \field[x_0,\ldots,x_N]$ be a set of homogeneous polynomials of degree at most $D$ and $J$ the homogeneous ideal generated by $F$. For a homogeneous ideal $I$ of $\field[x_0,\ldots, x_N]$  such that $\Z(I)$ is pure dimensional,  call $\mathcal{S}$ the set of irreducible components $Z$ of $\Z(I,J)$ with the property that $\Oc_{Z,\Z(I)}$ is Cohen-Macaulay. Then
 \begin{equation*}
 \sum_{Z\in \mathcal{S}}l_{Z,\Z(I)}(I,J)\cdot\deg(Z)\leq D^{\dim(\Z(I))-\dim (\Z(I,J))}\cdot\deg(\Z(I)).
 \end{equation*} 
 \end{thm}
 \begin{proof}
 See \cite[Ex. 12.3.7]{Fulton} or \cite[Thm. 1.1\text{ and } Remark 1.3]{Nakamaye1}.
 \end{proof}
We shall use the following consequence of Theorem \ref{R26}.
 \begin{cor}
 \label{R27}
  Let $F\subseteq \field[x_0,\ldots,x_N]$ be a set of homogeneous polynomials of degree at most  $D$ and $J$ the homogeneous ideal generated by $F$. Call $\mathcal{S}$ the set of irreducible components $Z$  of $\Z(I(\G),J)$ with the property that $Z\cap\phi(G)\neq\emptyset$. Then
 \begin{equation*}
 \sum_{Z\in \mathcal{S}}l_{Z,\phi(\G)}\big(I(\G),J\big)\cdot \deg(Z)\leq D^{n-\dim(\Z(I(\G),J))}\cdot \deg\big(\phi(\G)\big).
 \end{equation*} 
\end{cor}
\begin{proof}
From Theorem \ref{R26} it is enough to show that $\Oc_{Z,\phi(\G)}$ is Cohen-Macaulay if $Z\cap\phi(G)\neq\emptyset$. $G$ is smooth, see \cite[Sec. 1.2]{Borel}; in particular, $\Oc_{\phi(g),\phi(\G)}$ is Cohen-Macaulay for all $g\in G$. The localization of a Cohen-Macaulay ring by a prime ideal is Cohen-Macaulay, see \cite[Prop. 18.8]{Eisenbud}. Assume that $\phi(g)\in Z\cap\phi(G)$, then the local ring $\Oc_{Z,\phi(\G)}$ is isomorphic to the localization of $\Oc_{\phi(g),\phi(\G)}$ in the prime ideal corresponding to $Z$ and  then it  is Cohen-Macaulay by the previous argument.   
\end{proof}

Let $I$ be an ideal  of  $\field[x_0,\ldots, x_N]$ and denote by $\mathrm{Ass}\big(\field[x_0,\ldots,x_N]/I\big)$ the set of associated primes of the $\field[x_0,\ldots, x_N]-$module $\field[x_0,\ldots,x_N]/I$. The following result is a straight consequence of \cite[Lemma 3  and Prop. 1]{Wustholz1}; we emphasize that these results were stated for commutative algebraic groups, nevertheless their proofs work in the same way for noncommutative algebraic groups.\footnote{We just have to be careful in  \cite[ (8) p. 482]{Wustholz1}. The claim is true in the noncommutative case but the induction proof has to be done slightly more carefully.}
  \begin{lem}
 \label{R28}
 Let $I$ and $J$ be homogeneous ideals of $\field[x_0,\ldots,x_N]$ such that $J$ contains  $I(\G)$ and  $J\in\mathrm{Ass}\big(\field[x_0,\ldots,x_N]/I\big)$. If $\partial_{L_1}^TI\subseteq J$, then
 \begin{equation*}
 l_{\Z(J),\phi(\G)}(I)\geq \binom{\tau(\Z(J))+T}{\tau(\Z(J))}.
 \end{equation*}
 \end{lem}
 \begin{proof}
  See \cite[Lemma 3  and Prop. 1]{Wustholz1}. 
 \end{proof}
 
 Nakamaye \cite[Lemma 1.8]{Nakamaye1} gives a short proof the well known fact:  $\deg(\phi(V))=\deg(\phi(gV))$  for $V$ an irreducible variety of $\G$ and $g\in G$. In the same way, it can be proven that  $\deg(\phi(V))=\deg(\phi(Vg))$.
 \begin{lem}
 \label{R29}
 Let $V$ be an irreducible variety of $\G$. Then
 \begin{equation*}
 \deg(\phi(V))=\deg(\phi(gV))=\deg(\phi(Vg))\qquad\forall\; g\in G.
 \end{equation*}
 \end{lem}
 \begin{proof}
 See \cite[Lemma 1.8]{Nakamaye1}.
 \end{proof}
 
  Recall that given a  finite set $\Sigma_1$ of $G$ containing $1$ and $S\in\Nat$, $\Sigma_S$ is the set of products of $S$ elements of $\Sigma_1$. For $g\in \Sigma_S$ write $\phi(g)=\big[g_0:\ldots:g_N\big]$ and 
\begin{equation*}  
   J_{g}:=\Big(g_l\cdot x_k-g_k\cdot x_l:\;l,k\in\{0,\ldots,N\}\Big);
\end{equation*}  
  thus $J_g$ is the maximal ideal corresponding to $\phi(g)$ and $\Z(J_g)=\big\{\phi(g)\big\}$. Let $\Ac$ be the set of functions $f:\Sigma_S\rightarrow\{0,\ldots,N\}$  and write
  \begin{equation*}
  I_S:=\Bigg(\prod_{g\in \Sigma_S}(g_{e(g)}\cdot x_{f(g)}-g_{f(g)}\cdot x_{e(g)}):\;e,f\in \Ac\Bigg);
  \end{equation*}
     hence 
 \begin{equation}
 \label{E43}
 \Z(I_S)=\Z\Bigg(\prod_{g\in\Sigma_S} J_g\Bigg)=\Sigma_S.
 \end{equation}
 Define the ideals
\begin{align*}
\Id_{g,T}&:=\Big(P\in\field[x_0,\ldots, x_N]:\;P \text{ homogeneous},\; \ord_{g}(\bfr, P)\geq T\Big)\\
\Id_{S,T}&:=\Big(P\in\field[x_0,\ldots, x_N]:\;P \text{ homogeneous},\; \ord_{g}(\bfr, P)\geq T\quad\forall\;g\in \Sigma_S\Big).
\end{align*} 
 The next proposition is a trivial but fundamental tool to find obstruction subvarieties  in Theorem \ref{R1} and  Theorem \ref{R3}.
 \begin{prop}
 \label{R30}
 For all non-zero dimensional irreducible subvariety $W$ of $\G$, there is  $P\in\,^h\Id_{S,T}$ of degree at most $\sum_{i=0}^S(|\Sigma_1|-1)^i$ such that $\phi(W)\nsubseteq \Z(P)$.
 \end{prop}
 \begin{proof}
Let $\sqrt{I}$ denotes the radical of the ideal $I$. For all $g\in\Sigma_S$ the ideal $\Id_{g,1}$ is  prime so $\sqrt{\Id_{g,T}}\subseteq \sqrt{\Id_{g,1}}=\Id_{g,1}$. Note that $P(\phi(g))=0$ if  $P\in \Id_{g,1}$.  Then the Leibniz's rule implies that $\ord_{g}(\bfr, P^{T+1})\geq T$ if $P\in \Id_{g,1}$; this yields the inclusion $\Id_{g,1}\subseteq \sqrt{\Id_{g,T}}$ and we conclude that 
\begin{equation}
\label{E44}
\Id_{g,1}=\sqrt{\Id_{g,T}}.
\end{equation}
Hence 
\begin{align}
\label{E45}
\Z(\Id_{S,T})&=\Z\Bigg(\bigcap_{g\in\Sigma_S}\Id_{g,T}\Bigg)\nonumber\\
&=\bigcup_{g\in\Sigma_S}\Z(\Id_{g,T})\nonumber\\
&=\bigcup_{g\in \Sigma_S}\Z\big(\sqrt{\Id_{g,T}}\big)\nonumber\\
&=\bigcup_{g\in \Sigma_S}\Z(\Id_{g,1})&\text{ by (\ref{E44})}\nonumber\\
&=\Sigma_S\nonumber\\
&=\Z(I_S)&\text{by (\ref{E43}).}
\end{align}
Since $\dim(W)>0=\dim(\Z(I_S))$, we conclude from  (\ref{E45}) that there is one generator $P(x_0,\ldots,x_N):=\prod_{g\in \Sigma_S}(g_{e(g)}\cdot x_{f(g)}-g_{f(g)}\cdot x_{e(g)})$ of $I_S$ such that $\Z(P)\nsupseteq \phi(W)$. Since  $|\Sigma_S|\leq \sum_{i=0}^S(|\Sigma_1|-1)^i$, we conclude that $\deg(P)\leq \sum_{i=0}^S(|\Sigma_1|-1)^i$.
 \end{proof}
 Now we characterize $\ord_{g}(\bfr,T)$ in terms of $\partial^T_{L_g}(P)$  and $\partial^T_{R_g}(P)$. The following statement is analogous to \cite[Prop. 4.4]{Philippon}.
 \begin{cor}
 \label{R31}
 Let $P\in\field[x_0,\ldots,x_N]\setminus I(\G)$ be homogeneous, $g,h\in G$ and $T\in\Nat\cup\{0\}$. The following statements  are equivalent
 \begin{enumerate}
  \item[i)]$\ord_{gh}(\bfr,P)>T$.
 
  \item[ii)]$\phi(h)\in\Z\big(\partial^T_{L_g}(P)\big)$.
  
 \item[iii)]$\phi(g)\in\Z\big(\partial^T_{R_h}(P)\big)$.

 \end{enumerate}
 \end{cor}
 \begin{proof}
 Call $J:=\big(P_{\Delta,\alpha},I(\G):\;\alpha\in \Ac_R,\Delta\in\U(\bfr,T)\big)$. By the definition of the polynomials $P_{\Delta,\alpha}$, it is clear that $\ord_{gh}(\bfr,P)>T$ if and only if $gh\in\Z(J)$. Thus  the equivalence of i) and ii) is a consequence of the following equalities
  \begin{align*}
 \Z(J) &= \Z(\Eb^T(P))&\text{by  Corollary \ref{R10} iv)}\\
  &=\Z\big(\Tb_{L_{g^{-1}}}\big(\Ib\big(\partial^T_{L_g}(P)\big)\big)\big)&\text{by Proposition \ref{R22}}\\
  &=\phi\Big(g\cdot \phi^{-1}\big(\Z\big(\Ib\big(\partial^T_{L_g}(P)\big)\big)\big)\Big)&\text{by  Remark \ref{R12}}\\
  &=\phi\Big(g\cdot \phi^{-1}\big(\Z\big(\partial^T_{L_g}(P)\big)\big)\Big)&\text{by  Corollary \ref{R10} iv).}
  \end{align*}
 Likewise it is proven that i) and iii) are equivalent. 
 \end{proof}
 
 \section{Proofs of Theorem \ref{R1} and Theorem \ref{R2}}
 In this section we demonstrate Theorem \ref{R1} and Theorem \ref{R2}. Before we start with the proof of Theorem \ref{R1}, fix $P_1,\ldots,P_t \in\field[x_0,\ldots,x_N]$  homogeneous polynomials  such that $I(\G)=(P_1,\ldots,P_t)$ and define 
 \begin{align*}
 c_{7}:=\max_{1\leq k\leq t}\deg(P_k)\qquad\text{and}\qquad
 c_1=c_2=c_{5}^{2n}c_7^n\deg(\phi(\G)).
 \end{align*}
 \begin{proof}\emph{(Theorem \ref{R1})}
Let $\partial^T_{L_g}(I)$ and $\Id_{S,T}$ be as in Section 4 and Section 5. Denote by $I_2^*$ the ideal generated by $\bigcup_{g\in \Sigma_{[\frac{S}{n}]}}\partial^{[\frac{T}{n}]}_{L_{g}}(P)$. If $\dim(\Z(I_2^*))=n-1$, let $W_{2,1},\ldots, W_{2,m_2}$ be the irreducible components of  $\Z(I_2^*)$ of dimension $n-1$. From Proposition \ref{R30}, if $n-1>0$, there are homogeneous polynomials $Q_{2,1},\ldots, Q_{2,m_2}\in\Id_{S,T+1}$ of degree at most $\sum_{j=0}^S(|\Sigma_1|-1)^j$ such that $\phi(W_{2,i})\nsubseteq \Z(Q_{2,i})$ for each $ i\in\{1,\ldots, m_2\}$. Define
 \begin{equation*}
\mathcal{P}_2:= \left\{ \begin{array}{ll}
\emptyset & \mbox{if }\dim(\Z(I_2^*))<n-1,\;\dim(\Z(I_2^*))=0,\\
&\quad\text{ or }d_0>n-2;\\
&\\
\{Q_{2,1},\ldots, Q_{2,m_2}\}& \mbox{otherwise}.\end{array} \right.
\end{equation*}
We call $I_2$ the ideal generated by $I_2^*$ and $\mathcal{P}_2$.   Note that $\dim(\Z(I_2))<n-1$ or $\dim(\Z(I_2))=0$ insomuch as $\phi(W_i)\nsubseteq \Z(Q_{2,i})$ for all $ i\in\{1,\ldots m_2\}$.  We proceed with the construction of $\mathcal{P}_{r+1}$,  $I_{r+1}^*$ and $I_{r+1}$.   For $r\in\{2,\ldots,n\}$ let $I_{r+1}^*$ be the ideal generated by
 \begin{equation*}
\bigcup_{g\in\Sigma_{[\frac{rS}{n}]}}\partial^{[\frac{rT}{n}]}_{L_{g}}(P)\cup  \bigcup_{k=2}^r \,\bigcup_{Q\in\mathcal{P}_k}\,\bigcup_{g\in\Sigma_{[\frac{rS}{n}]}}\partial^{[\frac{rT}{n}]}_{L_{g}}(Q). 
 \end{equation*}
  If $\dim(\Z(I_{r+1}^*))=\dim(\Z(I_r))>0$, let $W_{r+1,1},\ldots, W_{r+1,m_{r+1}}$ be the irreducible components of  $\Z(I_{r+1}^*)$ of dimension $\dim(\Z(I_r))$. From Proposition \ref{R30} we know that there are homogeneous polynomials  $Q_{r+1,1}, \ldots, Q_{r+1,m_{r+1}} \in\Id_{S,T+1}$ of  degree  at most $\sum_{j=0}^S(|\Sigma_1|-1)^j$ such that $\phi(W_{r+1,i})\nsubseteq \Z(Q_{r+1,i})$ for each $i\in\{1,\ldots, m_{r+1}\}$. Define
 \begin{equation*}
\mathcal{P}_{r+1}:= \left\{ \begin{array}{ll}
\emptyset & \mbox{if }\dim(\Z(I_{r+1}^*))=0,\\
&\quad \dim(\Z(I_{r+1}^*))<\dim(\Z(I_r)),\\
&\quad\text{ or }d_0>n-r-1;\\
&\\
\{Q_{r+1,1},\ldots, Q_{r+1,m_{r+1}}\} & \mbox{otherwise}.\end{array} \right.
\end{equation*}
Call $I_{r+1}$  the ideal generated by $I_{r+1}^*$ and $\mathcal{P}_{r+1}$. The choice of the $Q_{r+1,i}$  let us conclude that if $r+1\leq n-d_0$, then   $\dim(\Z(I_{r+1}))<\dim(\Z(I_{r}))$ or $\dim(\Z(I_{r+1}))=0$; in particular $\dim(\Z(I_r))\leq n-r$ for all $r\leq n-d_0$. $I_{r+1}^*$  is generated by homogeneous polynomials of degree at most  $\max\big\{c_7,  c_5^2D\big\}$ and $\sum_{j=0}^S(|\Sigma_1|-1)^j\leq D$; then  $I_{r+1}$ is generated by  homogeneous  polynomials of degree at most $\max\big\{c_7,c_5^2D\big\}$. By Corollary \ref{R31},  $\phi(1)\in \Z\big(\partial^{T}_{L_g}(P)\big)$ and $\phi(1)\in \Z\Big(\bigcup_{Q\in\mathcal{P}_k}\partial^{T}_{L_{g}}(Q)\Big)$ for all $g\in\Sigma_S$ and $k\in\{2,\ldots, n\}$; this yields  $\phi(1)\in\Z(I_{n+1})\cap \phi(G)$. Let $d_r$  the maximal dimension of the irreducible components of $\Z(I_{r})$ which contain $\phi(1)$ for each $r\in\{2,\ldots,n+1\}$. Since
\begin{equation*}
\{\phi(1)\}\subseteq \Z(I_{n+1})\subseteq\ldots \subseteq \Z(I_{n-d_0+1})\subseteq \Z(I_{n-d_0})\subseteq\ldots \subseteq \Z(I_2)
\end{equation*}
and $d_r\leq\dim(\Z(I_r))\leq n-r$  when $r\leq n-d_0$, the Pigeonhole Principle yields that there is $r_0\in\{n-d_0,\ldots, n\}$ such that $d_{r_0}=d_{r_0+1}\leq d_0$; in particular, there is an irreducible  component  $W'$  of $\Z(I_{r_0})$ with dimension $d_{r_0}$ which is also an irreducible component of $\Z(I_{r_0+1})$ and  $\phi(1)\in W'\cap \phi(G)$. Let $J$ be the homogeneous prime ideal corresponding to $W'$ and write $W:=\phi^{-1}(W')$. The properties i), ii) and iii) are satisfied by the construction of $W$. It remains to prove that  iv) is true. The construction of $I_{r_0+1}$ and Corollary \ref{R25} i) lead to
\begin{equation}
\label{E46}
\partial^{[\frac{T}{n}]}_{L_g}(I_{r_0})\subseteq \Ib\Big(\partial^{[\frac{T}{n}]}_{L_g}(I_{r_0})\Big) \subseteq \Ib(I_{r_0+1})\subseteq J\qquad \forall\; g\in \Sigma_{[\frac{S}{n}]}.
\end{equation}
For all $g\in\Sigma_{[\frac{S}{n}]}$, let $J_g$ be the homogeneous prime ideal corresponding to the irreducible variety $\phi(gW)$. For all $g\in\Sigma_{[\frac{S}{n}]}$
\begin{align}
\label{E47}
\partial^{[\frac{T}{n}]}_{L_1}  (I_{r_0})&\subseteq\Ib\Big(\partial^{[\frac{T}{n}]}_{L_1}  (I_{r_0})\Big)\nonumber\\
&=\Ib\Big(\partial^0_{L_{g^{-1}}}\Big(\partial^{[\frac{T}{n}]}_{L_{g}}(I_{r_0})\Big)\Big)&\text{by Corollary \ref{R25}}\nonumber\\
&\subseteq \Ib\Big(\partial^0_{L_{g^{-1}}}(J)\Big)&\text{by (\ref{E46})}\nonumber\\
&\subseteq J_g&\text{by Remark \ref{R12}}.
\end{align}
Then Lemma \ref{R28} and (\ref{E47}) imply that
\begin{equation}
\label{E48}
\binom{\tau(gW)+\big[\frac{T}{n}\big]}{\tau(gW)}\leq l_{\phi(gW), \phi(\G)}(I_{r_0}) \qquad\forall\;g\in\Sigma_{[\frac{S}{n}]}.
\end{equation} 
Since the left translations are isomorphisms
\begin{align}
\label{E49}
\tau(W)=\tau(gW).
\end{align}
  Let $S_W$ be the set of different irreducible varieties in $\Big\{gW:\;g\in\Sigma_{[\frac{S}{n}]}\Big\}$ so $|S_W|=N_W$. To conclude the proof of the theorem,  name $r_W:=N_W\binom{[\frac{T}{n}]+\tau(W)}{\tau(W)}\deg(\phi(W))$ and see
\begin{align*}
r_W&=\sum_{gW\in S_W}\binom{\big[\frac{T}{n}\big]+\tau(gW)}{\tau(gW)}\deg(\phi(W))&\text{by (\ref{E49})}\\
&\leq  \sum_{gW\in S_W}l_{\phi(gW), \phi(\G)}(I_{r_0})\deg(\phi(W))&\text{by (\ref{E48})}\\
&= \sum_{gW\in S_W} l_{\phi(gW), \phi(\G)}(I_{r_0})\deg(\phi(gW))&\text{by Lemma \ref{R29}}\\
&\leq \max\{c_7,c_5^2D\}^{n-\dim(\Z(I_{r_0}))} \deg(\phi(\G))&\text{by Corollary \ref{R27}}\\
&\leq c_1D^{n-\dim(W)}. 
\end{align*}
 \end{proof}
 Now we proceed with the proof Theorem \ref{R2}. Broadly speaking, the main idea of the proof of Theorem \ref{R2} is the same idea used Theorem \ref{R1}; furthermore, since we do not have the assumption $\sum_{i=0}^S(|\Sigma_1|-1)^i\leq D$, we wont find an upper bound of the dimension of the obstruction variety and  this fact makes the proof of Theorem \ref{R2} easier.
 \begin{proof}\emph{(Theorem \ref{R2})}
 For all $r\in\{2,\ldots, n+1\}$ denote by $I_r$ the homogeneous ideal generated by $\bigcup_{g\in\Sigma_{[\frac{(r-1)S}{n}]}}\partial^{[\frac{(r-1)T}{n}]}_{L_{g}}(P)$.    $I_{r}$  is generated by polynomials of degree at most  $\max\big\{c_{7},  c_5^2D\big\}$. Call $d_r$ the maximal dimension of the irreducible components of $\Z(I_r)$ containing $1$. By Corollary \ref{R31}, $\phi(1)\in \Z\big(\partial^{T}_{L_{g}}(P)\big)$  for all $g\in\Sigma_S$ and  this gives that $\phi(1)\in\Z(I_{n+1})\cap \phi(G)$. Moreover, the inclusions
\begin{equation*}
\{\phi(1)\}\subseteq \Z(I_{n+1})\subseteq\ldots \subseteq \Z(I_2)
\end{equation*}
and the Pigeonhole Principle let us conclude that there is $r_0\in\{2,\ldots, n\}$ such that $d_{r_0}=d_{r_0+1}$; in particular, there is an irreducible  component  $W'$  of $\Z(I_{r_0})$ which is also an irreducible component of $\Z(I_{r_0+1})$ and  $\phi(1)\in W'\cap \phi(G)$. Name $W:=\phi^{-1}(W')$ and note that the properties of $W'$ yields i) and ii) are satisfied. The fact that $W$ satisfies iii) is proven exactly in the same way as the last part of the proof of Theorem \ref{R1}.
 \end{proof}
  \section{Proofs of Theorem \ref{R3} and Theorem \ref{R4}}
 We demonstrate Theorem \ref{R3} and Theorem \ref{R4} in this section. The main idea that we will follow is quite similar to the ones taken in the proofs of Theorem \ref{R1} and Theorem \ref{R2}; nonetheless, to show that the obstruction variety that we obtain is the closure of an algebraic group, it is not enough to conclude as we did in the proof of Theorem \ref{R1}. We start with the proof of Theorem \ref{R3} where we will see that it suffices to take $c_3=c_{5}^{3n}c_7^n\deg(\phi(\G))$; in the first part of the proof, we use a similar construction as the one done in Theorem \ref{R1}; in the second part we conclude with the same ideas of the last part of  \cite[Lemme 5.1]{Philippon} and \cite[Thm. 0.3]{Nakamaye1}. 
\begin{proof}\emph{(Theorem \ref{R3})}
Since $gB=Bg$ for all $g\in \Sigma_1$, we have that $\mathrm{Ad}\big(g\big)(\bfr)\subseteq \bfr$ for all $g\in\Sigma_S$, see \cite[Sec. 9.2]{Hilgert-Neeb}. Denote by $I_2^*$ the ideal generated by $\bigcup_{g\in\Sigma_{[\frac{S}{n}]}}\partial^{[\frac{T}{n}]}_{R_{g}}(P)$. If $\dim(\Z(I_2^*))=n-1$, let $W_{2,1},\ldots, W_{2,m_2}$ be the irreducible components of  $\Z(I_2^*)$ with dimension $n-1$. Proposition \ref{R30} assures the existence, if $n-1>0$, of homogeneous polynomials $Q_{2,1},\ldots, Q_{2,m_2}\in\Id_{S,T+1}$ of degree at most $\sum_{j=0}^S(|\Sigma_1|-1)^j\leq D$ such that $\phi(W_{2,i})\nsubseteq \Z(Q_{2,i})$ for each $i\in\{1,\ldots, m_2\}$. Define
 \begin{equation*}
\mathcal{P}_2:= \left\{ \begin{array}{ll}
\emptyset & \mbox{if }\dim(\Z(I_2^*))<n-1,\;\dim(\Z(I_2^*))=0,\\
&\quad\text{ or }d_0>n-2;\\
&\\
\{Q_{2,1},\ldots, Q_{2,k_2}\}& \mbox{otherwise}.\end{array} \right.
\end{equation*}
We call $I_2$ the ideal generated by $I_2^*$ and $\mathcal{P}_2$. Thus $\dim\Z(I_2)<n-1$ or $\dim\Z(I_2)=0$. Now take $r\in\{2,\ldots,n\}$ and  let $I_{r+1}^*$ be the ideal generated by
 \begin{equation*}
\bigcup_{g\in\Sigma_{[\frac{rS}{n}]}}\partial^{[\frac{rT}{n}]}_{R_{g}}(P)\cup  \bigcup_{k=2}^r\, \bigcup_{Q\in\mathcal{P}_k}\,\bigcup_{g\in\Sigma_{[\frac{rS}{n}]}}\partial^{[\frac{rT}{n}]}_{R_{g}}(Q). 
 \end{equation*}
  If $\dim(\Z(I_{r+1}^*))=\dim(\Z(I_r))>0$, let $W_{r+1,1},\ldots, W_{r+1,m_{r+1}}$ be the irreducible components of  $\Z(I_{r+1}^*)$ of dimension $\dim(\Z(I_r))$. From Proposition \ref{R30} there are homogeneous polynomials  $Q_{r+1,1}, \ldots, Q_{r+1,m_{r+1}} \in\Id_{S,T+1}$ of  degree  at most $D$ such that $\phi(W_{r+1,i})\nsubseteq \Z(Q_{r+1,i})$ for each $i\in\{1,\ldots, m_{r+1}\}$. Define
 \begin{equation*}
\mathcal{P}_{r+1}:= \left\{ \begin{array}{ll}
\emptyset & \mbox{if }\dim(\Z(I_{r+1}^*))=0,\\
&\quad \dim(\Z(I_{r+1}^*))<\dim(\Z(I_r)),\\
&\quad\text{ or }d_0>n-r-1;\\
&\\
\{Q_{r+1,1},\ldots, Q_{r+1,k_{r+1}}\} & \mbox{otherwise}.\end{array} \right.
\end{equation*}
Call $I_{r+1}$  the ideal generated by $I_{r+1}^*$ and $\mathcal{P}_{r+1}$. The construction of $\mathcal{P}_{r+1}$ implies that if $r+1\leq n-d_0$, then $\dim\Z(I_{r+1})<\dim\Z(I_r)$ or $\dim\Z(I_{r+1})$. As a consequence of Corollary \ref{R31}, $\phi(1)\in \Z\big(\partial^{T}_{R_{g}}(P)\big)$ and $\phi(1)\in \Z\Big(\bigcup_{Q\in\mathcal{P}_k}\partial^{T}_{R_{g}}(Q)\Big)$ for all $g\in\Sigma_S$ and $k\in\{2,\ldots, n\}$; thus $\phi(1)\in\Z(I_{n+1})\cap \phi(G)$. Let $d_r$ the maximal dimension of the irreducible components of $\Z(I_{r})$ which contain $\phi(1)$. Since
\begin{equation*}
\{\phi(1)\}\subseteq \Z(I_{n+1})\subseteq\ldots \subseteq \Z(I_{n-d_0+1})\subseteq \Z(I_{n-d_0})\subseteq\ldots \subseteq \Z(I_2)
\end{equation*}
and $d_r\leq\dim(\Z(I_r))\leq n-r$ when $r\leq n-d_0$, the Pigeonhole Principle yields the existence of $r_0\in\{n-d_0,\ldots, n\}$ such that $d_{r_0}=d_{r_0+1}\leq d_0$; in particular, there is an irreducible  component  $V'$  of $\Z(I_{r_0})$ which is also an irreducible component of $\Z(I_{r_0+1})$ with $\dim V'\leq d_0$ and  $\phi(1)\in V'\cap \phi(G)$. Set $V:=\phi^{-1}(V')$ and $W':=\bigcap_{g\in G}gVg^{-1}$. Let $W$ be the irreducible component of $W'$ which contains $1$ and hence $gWg^{-1}=W$ for all $g\in G$. Now define $H':=\{h\in G:\; Wh=W\}$ and note that $H'$ is an algebraic subgroup since $W\cap G$ is  closed in $G$, see \cite[Sec. 1.1]{Borel}. Insomuch as $1\in W$
\begin{equation}
\label{E50}
\dim(H')\leq \dim(W)\leq d_0. 
\end{equation}
Call $I$ the homogeneous ideal of $\field[x_0,\ldots,x_N]$ generated by $\bigcup_{h\in W\cap G}\partial^0_{L_h}(I_{r_0})$; particularly, $I$ is generated by homogeneous polynomials of degree at most $\max\big\{c_7,c_5^3D\big\}$. Denote by $H$ the irreducible component of $H'$ which contains $1$. If $S$ is a subset of $G$, $\overline{S}$ is the closure of $S$ in $\G$;   then  $\phi(\overline{H})\subseteq \Z(I)$. From (\ref{E50}) and the trivial inclusion $\Z(I)\subseteq \Z(P)$, the conditions i) and ii) are proven. See that $H$ is a normal subgroup of $G$ since $gWg^{-1}=W$ for all $g\in G$. Then it remains to show  iv). For all $g\in\Sigma_{[\frac{S}{n}]}$ let $J_g$ be the homogeneous prime ideal  corresponding to $\phi\big(\overline{Hg}\big)=\phi\big(\overline{H}\cdot g\big)$.  The construction of $I_{r_0+1}$ leads to the inclusion
\begin{equation}
\label{E51}
 \partial^{[\frac{T}{n}]}_{R_{g}}(I_{r_0})\subseteq I_{r_0+1}\qquad \forall\;g\in\Sigma_{[\frac{S}{n}]}.
\end{equation}
Since $WH=W$, we get that 
\begin{equation*}
(W\cap G)\cdot\overline{H}\subseteq W
\end{equation*}
and consequently
\begin{equation}
\label{E52}
\Ib\big(\partial^0_{L_h}(I_{r_0+1})\big)\subseteq J_1\qquad \forall\;h\in W\cap G.
\end{equation}
Then
\begin{align}
\label{E53}
\partial^{[\frac{T}{n}]}_{R_{g}}(I)&=\partial^{[\frac{T}{n}]}_{R_{g}}\Bigg(\bigcup_{h\in W\cap G}\partial^0_{L_h}(I_{r_0})\Bigg)\nonumber\\
&\subseteq\Ib\Bigg(\bigcup_{h\in W\cap G}\partial^{[\frac{T}{n}]}_{R_{g}}\big(\partial^0_{L_h}(I_{r_0})\big)\Bigg)\nonumber\\
&=\Bigg(\bigcup_{h\in W\cap G}\Ib\Big(\partial^{[\frac{T}{n}]}_{R_{g}}\big(\partial^0_{L_h}(I_{r_0})\big)\Big)\Bigg)&\text{by Corollary \ref{R10} ii)}\nonumber\\
&=\Bigg(\bigcup_{h\in W\cap G}\Ib\Big(\partial^0_{L_h}\Big(\partial^{[\frac{T}{n}]}_{R_{g}}(I_{r_0})\Big)\Big)\Bigg)&\text{by Corollary \ref{R25}}\nonumber\\
&\subseteq \Bigg(\bigcup_{h\in W\cap G}\Ib\big(\partial^0_{L_h}(I_{r_0+1})\big)\Bigg)&\text{by (\ref{E51})}\nonumber\\
&\subseteq J_1&\text{by (\ref{E52}).}
\end{align}
and for all $g\in\Sigma_{[\frac{S}{n}]}$  
\begin{align}
\label{E54}
\partial^{[\frac{T}{n}]}_{L_1} (I)&\subseteq \Ib\Big(\partial^{[\frac{T}{n}]}_{L_1} (I)\Big)\nonumber\\
&=\Ib\Big(\partial^{[\frac{T}{n}]}_{R_1} (I)\Big)&\text{by Lemma \ref{R14} iv)}\nonumber\\
&=\Ib\Big(\partial^0_{R_{g^{-1}}}\Big(\partial^{[\frac{T}{n}]}_{R_{g}}(I)\Big)\Big)&\text{by Corollary \ref{R25}}\nonumber\\
&=\Ib\Big(\partial^0_{R_{g^{-1}}}(J_1)\Big)&\text{by (\ref{E53})}\nonumber\\
&\subseteq J_{g}. 
\end{align}
 Lemma \ref{R28} and (\ref{E54}) imply that
\begin{equation}
\label{E55}
\binom{\tau\big(\overline{Hg}\big)+\big[\frac{T}{n}\big]}{\tau\big(\overline{Hg}\big)}\leq l_{\phi(\overline{Hg}), \phi(\G)}(I) \qquad\forall\;g\in\Sigma_{[\frac{S}{n}]}.
\end{equation} 
Also for all $g\in \Sigma_{[\frac{S}{n}]}$
\begin{align}
\label{E56}
\tau\big(\overline{H}\big)=\tau\big(\overline{Hg}\big).
\end{align}
  Let $S_H$ be the set of different cosets in $\big\{Hg:\;g\in\Sigma_{[\frac{S}{n}]}\big\}$ so $|S_H|=N_H$. Call $r_H:=N_H\binom{[\frac{T}{n}]+\tau(\overline{H})}{\tau(\overline{H})}\deg(\phi(\overline{H}))$. Then
\begin{align*}
r_H&=\sum_{Hg\in S_H}\binom{\big[\frac{T}{n}\big]+\tau\big(\overline{Hg}\big)}{\tau\big(\overline{Hg}\big)}\deg(\phi(\overline{H}))&\text{by (\ref{E56})}\\
&\leq  \sum_{Hg\in S_H}l_{\phi(\overline{Hg}), \phi(\G)}(I)\deg(\phi(\overline{H}))&\text{by (\ref{E55})}\\
&\leq \sum_{Hg\in S_H} l_{\phi(\overline{Hg}), \phi(\G)}(I)\deg(\phi(\overline{Hg}))&\text{by Lemma \ref{R29}}\\
&\leq \max\big\{c_7,c_5^3D\big\}^{n-\dim(\Z(I))} \deg(\phi(\G))&\text{by Corollary \ref{R27}}\\
&\leq c_3D^{n-\dim(H)}. 
\end{align*}
\end{proof}
Now we  show Theorem \ref{R4} with $c_4=c_{5}^{3n}c_7^n\deg(\phi(\G))$.  
 \begin{proof}\emph{(Theorem \ref{R4})}
 For all $r\in\{2,\ldots, n+1\}$, denote by $I_r$ the homogeneous ideal generated by $\bigcup_{g\in\Sigma_{[\frac{(r-1)S}{n}]}}\partial^{[\frac{(r-1)T}{n}]}_{R_{g}}(P)$.  Call $d_r$ the maximal dimension of the irreducible components of $\Z(I_r)$ containing $1$. From Corollary \ref{R31}, $\phi(1)\in \Z\big(\partial^{T}_{R_{g}}(P)\big)$ for all $g\in\Sigma_S$; this yields  $\phi(1)\in\Z(I_{n+1})\cap \phi(G)$. Furthermore, the inclusions
\begin{equation*}
\{\phi(1)\}\subseteq \Z(I_{n+1})\subseteq\ldots \subseteq \Z(I_2)
\end{equation*}
and the Pigeonhole Principle let us conclude the existence of $2\leq r_0\leq n$ such that $d_{r_0}=d_{r_0+1}$; in particular, there is an irreducible  component  $V'$  of $\Z(I_{r_0})$ which is also an irreducible component of $\Z(I_{r_0+1})$ with $\phi(1)\in V'\cap \phi(G)$ and write $V:=\phi^{-1}(V)$. Set $W':=\bigcap_{g\in G}gVg^{-1}$ and $W$ its maximal dimensional irreducible component containing $1$.  Call $H':=\{h\in G:\; Wh=W\}$ and $H$ the irreducible component of $H'$ containing $1$. The conclusion of the proof is exactly the same as the one of Theorem \ref{R3}.
 \end{proof}

\end{document}